\documentclass[dvips,ss]{imsart}

\RequirePackage[OT1]{fontenc}
\RequirePackage{amsthm,amsmath}
\RequirePackage[numbers]{natbib}
\RequirePackage[colorlinks,citecolor=blue,urlcolor=blue]{hyperref}
\usepackage{graphicx}
\pubyear{2012}
\volume{0} \issue{0}
\firstpage{1}
\lastpage{8}

\startlocaldefs
\numberwithin{equation}{section}
\theoremstyle{plain}

\newcommand{\df}{\stackrel{\bigtriangleup}{=}}
\endlocaldefs

\begin{document}

\begin{frontmatter}
\title{Multiple Objects: Error Exponents in  Hypotheses  Testing and    Identification}
\runtitle {Multiple Objects:   Hypotheses  Testing and
Identification}

\begin{aug}
\author{\fnms{Evgueni} \snm{Haroutunian}
\ead[label=e1]{evhar@ipia.sci.am}}\, \and\,
\author{\fnms{Parandzem} \snm{Hakobyan} \ead[label=e2]{par\_h@ipia.sci.am}}

\address{Institute for Informatics and Automation Problems\\
 of the National Academy of Sciences of the Republic of Armenia }

 \runauthor{E. Haroutunian et al.}

\end{aug}

\begin{abstract}
We  survey  a series of investigations of optimal testing of
multiple hypotheses concerning various multiobject models.

These studies are a bright instance of application of methods and
technics developed in Shannon information theory to solution of
typical statistical problems.
\end{abstract}


\end{frontmatter}

\section {Introduction}

 ``One can conceive of Information Theory in the broad sence as covering the theory of Gaining, Transferring, and
  Storing Information, where the first is usually called
  Statistics."\cite {Rnew}.

Shannon information theory and mathematical statistics interaction
revealed to be  effective. This interplay  is mutually fruitful, in
some works results of probability theory and statistics were
obtained with application of information-theoretical methods and
there are studies where statistical results provide ground for new
findings in information theory \cite {14}, \cite {16}, \cite
{3C}--\cite {17}, \cite {SE}, \cite {My6}, \cite {18}, \cite {10},
\cite {6}--\cite {19}.

This paper can serve an illustration of application of
information-theoretical methods in statistics: on  one hand this is
analogy in problem formulation and on the other hand this is
employment of technical tools of proof, specifically of the method
of types \cite {3}, \cite {2}.

 It is often necessary in statistical research to make decisions regarding the nature and parameters
 of stochastic model, in particular, the    probability
 distribution of the object.
 Decisions can be  made on the base of results of observations of the object.
 The vector of results is called a  sample.
 The  correspondence between samples and hypotheses can be  designed  based on some selected criterion.
 The procedure of statistical hypotheses detection is called test.

The classical  problem of statistical hypothesis testing  refers to
two hypotheses. Based on data samples a statistician makes decision
on which of the two proposed hypotheses must be  accepted. Many
mathematical investigations, some of which have also applied
significance, were implemented in this direction \cite {LEM}.

The need of  testing of more than two hypotheses in many scientific
and applied fields has essentially increased recently.  As an
instance microarray analysis could be mentioned \cite {28*}.

 The decisions can  be erroneous due to randomness of the sample.
 The test is considered as  good if the probabilities of the errors
 in given conditions are as small as possible.

Frequently the problem was solved for the case of a tests sequence,
where the probabilities of error decrease exponentially as
$2^{-NE}$, when the number of observations $N$ tends to the
infinity. We call the exponent of error probability $E$ the {\it
reliability}. In case of two hypotheses both reliabilities
corresponding to two possible error probabilities could not be
increased  simultaneously, it is an accepted way to fix the value of
one of the reliabilities and try to make the tests sequence get the
greatest value of the remaining reliability. Such a test is called
{\it logarithmically asymptotically optimal} (LAO). Such optimal
tests were considered first by Hoeffding \cite {4}, examined later
by Csisz\'ar and  Longo \cite {5}, Tusnady \cite {6}, \cite {7}
 (he called such test series {\it exponentially rate optimal} (ERO)),
  Longo and  Sgarro \cite {8}.
  The term  LAO  for  testing of two hypotheses was  proposed by Birge  \cite {9}. Amongst
papers   on testing,  associated with information theory,  we  can
also note works of Natarajan \cite {10}, Gutman \cite {11},
 Anantharam \cite {12}, Han  \cite {13} and of many others. Some  objectives in this direction were first
suggested in  original introductory article by Dobrushin, Pinsker
and Shiryaev \cite {20}. The  achievable region  of  error
probability exponents was examined by Tuncel  \cite {Tn}.

 The problem has common
features with the issue studied in the information theory on
interrelation between the rate $R$ of the code  and the exponent $E$
of  the error probability. In information theory the relation $E(R)$
is called according to Shannon the {\it reliability function}, while
$R(E)$ is named  the $E$-{\it capacity}, or the {\it
reliability-rate function}, as it was introduced  by Haroutunian
\cite {1}, \cite {H07}, \cite {HMA}.


Simple but actual concept of not only separate but also simultaneous
investigation of some number  of objects of the same type,
evidently, was first formulated by Ahlswede and Haroutunian \cite
{R} for reliable testing of distributions of multiple items.  But
simultaneous examination of properties of many similar objects may
be attractive and effective  in plenty of other statistical
situations.

 The organization of this paper is as
follows. We start with the definitions and notations in the next
section. In section 3 we introduce the problem of multihypotheses
testing concerning one object. In section 4 we consider the
reliability approach to multihypotheses testing for many independent
and dependent objects. Section 5 is dedicated to the problem of
statistical identification under condition of optimality. Section 6
is devoted to description of characteristics of LAO hypotheses
testing with permission of rejection of decision for the model
consisting of one and of  more independent objects.

\section {Definitions and Notations}

 We  denote finite sets by script capitals. The cardinality of a set
${\cal X}$ is denoted as  $|{\cal X}|$. Random variables (RVs),
which take values in finite sets ${\cal X}$, ${\cal S}$ are denoted
by $X$, $S$. Probability distributions (PDs) are denoted by $Q$,
$P$, $G$, $W$, $V$, $Q\,\mbox {o}  \,V$.

 Let PD of RV $X$, characterizing an object,  be $Q\df \{Q(x), \,\,x\in {\cal X}\}$, and
 conditional PD of RV $X$ for given value of state  $s$ of the object  be  $V\df
 \{V(x|s), \,x\in {\cal X}, \,s\in {\cal S}\}$.

 The Shannon  entropy $H_Q(X)$ of RV $X$   with PD $Q$ is:
\begin {eqnarray*}
H_Q(X)\df -\sum\limits_{x\in {\cal X}}Q(x)\log Q(x).
\end {eqnarray*}
 The conditional entropy $H_{P, V}(X\mid S)$ of RV $X$
for given RV $S$ with corresponding PDs  is:
\begin {align*}
H_{P, V}(X\mid S)\df -\sum\limits_{x\in {\cal X}, s\in {\cal
S}}P(s)V(x|s)\log {V(x| s)}.
\end {align*}

The  divergence (Kullback-Leibler information, or ``distance") of
PDs $Q$ and $G$ on ${\cal X}$ is:
$$
D(Q||G)\df \sum\limits_{x\in {\cal X}} Q(x)\log \frac {Q(x)}{G(x)},
$$
and   conditional divergence of the PD  $P\,\mbox {o}
 \,V=\!\!\{P(s)V(x|s),\,\, x\in {\cal X},\,\,s\in {\cal S}
\}$ and  PD $P\,\mbox {o}  \,W \!=\!\{P(s)\,W(x|s),x\in {\cal
X},\,\, s\in {\cal S} \}$ is:
$$
D(P\circ V||P\circ W )=D(V||W|P)\df
\sum\limits_{x,s}P(s)V(x|s)\log\frac{V(x|s)}{Wx|s)}.
$$

For our investigations we use the method of types,   one of the
important technical tools in Shannon theory \cite {2,3}. The type
$Q_{\bf x}$ of a vector ${\bf x}=(x_1,...,x_N)\in {\cal X}^N$ is a
PD (the empirical distribution)
$$
Q_{\bf x}=\left\{Q_{\bf x}(x)= \frac {N(x|{\bf x})}{N}, x\in {\cal
X}\right \},
$$
 where $N(x|{\bf x})$ is the number of repetitions of symbol $x$ in vector ${\bf x}$.

  The joint type of vectors ${\bf x}\in {\cal X}$ and  ${\bf s}=(s_1,s_2,...,s_N)\in
{\cal S}^N$ is the PD
$$P_{\bf x,\, s}=\left \{\frac {N(x, s|{\bf x}, {\bf s})}{N}, \,\,x\in {\cal
X},\,\,{  s}\in {\cal S}\right \},$$
 where $N(x, s|{\bf x}, {\bf s})$ is the number of occurrences of symbols pair $(x,s)$ in the pair of vectors
$({\bf x}, {\bf s})$.
 The conditional type of $\bf x$ for given $\bf
s$ is    conditional PD
 $$
 V_{{\bf x}|{\bf s}}=\{V(x|s), x\in {\cal X}, s\in {\cal S}\},
 $$
   defined by relation $N(x, s|{\bf x}, {\bf
s})=N(x |{\bf x})V_{{\bf x}|{\bf s}}(x|s)$ for all $x\in {\cal X}$,
${ s}\in {\cal S}.$

We denote by ${\cal Q}^N({\cal X})$ the set of all  types of vectors
in  ${\cal X}^N$ for given $N$,   by  ${\cal P}^N({\cal S})$ -- the
set of all types of vectors ${\bf s}$ in  ${\cal S}^N$
 and by ${\cal V}^N({\cal X}|{\bf s})$ -- the
set of all possible conditional types of vectors ${\bf x}$ in ${\cal
X}^N$ for given ${\bf s}\in {\cal S}^N$. The set of vectors $\bf x$
of type $Q$ is denoted by ${\cal T}_{Q}^N(X)$ and     the family of
vectors ${\bf x}$ of  conditional type $V$ for given ${\bf s}\in
{\cal S}^N$ of  type $P$ by ${\cal T}_{P,V}^{N}(X\mid {\bf s})$. The
set of all possible PDs $Q$ on  ${\cal X}$ and PDs $P$ on ${\cal S}$
is  denoted, correspondingly, by ${\cal Q}(X)$  and ${\cal P}({\cal
S} )$.

 We need the following frequently used  inequalities \cite {2}:
\begin {equation}
\label{1}\mid {\cal Q}^N({\cal X})\mid \leq (N+1)^{ \mid {\cal
X}\mid }, \end {equation}
\begin {equation}
\label{1*}\mid {\cal V}^N({\cal X}| {\bf s})\mid \leq (N+1)^{\mid
{\cal S}\mid \mid {\cal X}\mid },
\end {equation}
for any type $Q\in {\cal Q}^N({\cal X})$
\begin {equation}
\label {2} (N+1)^{-|{\cal X}|}\exp\{NH_{Q}(X)\}\leq\mid {\cal
T}_{Q}^{N}(X) \mid \leq\exp\{H_{Q}(X)\},
\end {equation}
 and  for any type
$P\in {\cal P}^N({\cal S})$ and $V \in {\cal V}^N({\cal X}|{\bf s})$
\begin {equation}
\label {2*} (N+1)^{-|{\cal S}||{\cal
X}|}\exp\{NH_{P,V}(X|S)\}\leq\mid {\cal T}_{P,V}^{N}(X\mid {\bf s})
\mid \leq\exp\{H_{P,V}(X|S)\}.
\end {equation}

\section {LAO Testing of Multiple    Hypotheses for One Object}

The problem of optimal testing of {\it multiple} hypotheses  was
proposed by Dobrushin \cite {D}, and  was investigated in  \cite
{E4}--\cite {E5}.  Fu and
  Shen \cite {FS} explored the case of two hypotheses  when side information is
  absent.
The problem   concerning arbitrarily varying sources solved in \cite
{My1} was induced by the ideas of the paper of Ahlswede \cite {R1}.
The case of two hypotheses with  side information about states was
considered in \cite {el}.
  In the same way as in \cite {FS}  from result on  LAO testing,   the
  rate-reliability and the reliability-rate functions for
  arbitrarily varying source with side information were obtained in
  \cite {My1}.

    The problem of multiple hypotheses LAO testing for discrete
stationary Markov source of observations  was solved by Haroutunian
\cite {Em1}--\cite {Em3}. In \cite {EN1} Haroutunian and Grigoryan
generalized results  from   \cite {FS}, \cite {Em1}--\cite {Em3} for
multihypotheses LAO testing by a non-informed statistician for
arbitrarily varying Markov source.

Here for clearness we expose the results on multiple hypotheses LAO
testing for the case of the  most simple invariant object.

   Let  ${\cal X}$ be a finite set of values of random variable (RV) $X$.
   $M$ possible PDs $ G_m=\{G_m(x), \,\,\,x\in {\cal X}\}$, $m=\overline
{1,M}$, of RV $X$ characterizing  the object are known.

 The statistician must detect  one among  $M$ alternative
hypotheses  $G_m$,  using     sample ${\bf x}=(x_1,...,x_N)$ of
results of $N$ independent observations of the object.

 The procedure of decision making is a non-randomized test $\varphi_N({\bf x})$,
 it can be defined  by division of the sample space ${\cal X}^N$ on $M$ disjoint subsets
  ${\cal A}^N_m=\{{\bf x:}\,\,\,
  \varphi_N({\bf x})=m\}$,   $m=\overline {1,M}$.
  The  set  ${\cal A}^N_m$ consists of all samples $\bf x$ for which the hypothesis
   $G_m$ must be  adopted.
 We study      the probabilities
 $\alpha_{l|m}(\varphi_N)$ of the erroneous acceptance
  of hypothesis $G_l$ provided that $G_m$ is true
\begin {equation}
\label {4} \alpha_{l|m}(\varphi_N)\df G^N_m(A_l^{N}),\,\,
l,m=\overline{1,M},\,\,\,m\not= l. \end {equation}

 The probability to reject hypothesis $G_m$, when it is true, is also considered
\begin {align}
 \alpha_{m|m}(\varphi_N)&  \df \sum\limits_{l\not
 =m}\alpha_{l|m}(\varphi_N)\nonumber\\
 & = G^N_m(\overline {{\cal A}^N_m})\nonumber\\
 & =(1-G^N_m({{\cal A}^N_m})).\label {5}
\end {align}

A quadratic matrix  of $M^2$ error probabilities
$\{\alpha_{l|m}(\varphi_N),\,\,m=\overline{1,M},\,\,l=\overline{1,M}\}$
  is  the {\it  power} of the  tests.

 {\it Error probability exponents} of the infinite
sequence $\varphi$ of tests, which we call  {\it reliabilities}, are
defined as follows:
\begin {equation}
\label {6} E_{l|m}(\varphi)\df\overline {\lim\limits_{N\to
\infty}}\left \{ -\frac {1}{N}\log\alpha_{l|m}(\varphi_N)\right
\},\,\,\,m,l=\overline {1,M}.
\end {equation}
 We see from (\ref {5}) and (\ref {6}) that
\begin {equation}
\label {7} E_{m|m}(\varphi)=\min_{l\not=
m}E_{l|m}(\varphi),\,\,\,m=\overline {1,M}.
\end {equation}
 The
matrix
$$
{\bf E}(\varphi)=\left (
\begin{array}{c}
E_{1| 1}\,\, \ldots \, E_{ l|1}\, \ldots \,\, E_{M|1 }\\
\ldots \ldots \ldots \ldots  \ldots \ldots \ldots \\
E_{ 1|m} \ldots E_{l|m} \ldots   E_{ M|m},\\
\ldots \ldots \ldots \ldots  \ldots \ldots \ldots \\
E_{1|M } \ldots  E_{ l|M } \ldots  E_{M| M}
\end{array}
\right )
$$
called the {\it reliabilities matrix} of the tests sequence
$\varphi$  is the object of our investigation.

   We recognize that a sequence $\varphi^*$ of tests is
LAO if for given positive values of $M-1$  diagonal, elements of
matrix ${\bf E}(\varphi^*)$ the procedure  provides maximal values
for all  other elements of it.

\vspace {2mm}
 Now we form   the LAO test by constructing   decision sets noted ${\cal R}_m^{(N)}$.
 Given strictly positive numbers $E_{m|m}$, $m=\overline {1,M-1}$, we define  the following regions:
\begin {equation}
\label {8} {\cal R}_{m}\df \{Q:\,\,\,D(Q||G_m)\leq
E_{m|m}\},\,\,\,\,m=\overline{1,M-1},
\end {equation}
\begin {equation}
\label {9} {\cal R}_{M} \df \{Q:\,\,\,D(Q||G_m)>
E_{m|m},\,\,\,\,m=\overline{1,M-1}\},\end {equation}
\begin {equation}
\label {10} {\cal R}^{(N)}_{m}\df {\cal R}_{m}\bigcap {\cal
Q}^N({\cal X}),\,\,\,m=\overline {1,M},
\end {equation}
and corresponding values:
\begin {equation}
\label {11} E^*_{m|m}=E^*_{m|m}(E_{m|m})\df
E_{m|m},\,\,\,\,\,m=\overline {1,M-1},
\end {equation}
\begin {equation}
\label {12} E^*_{m|l}=E^*_{m|l}(E_{m|m})\df \inf\limits_{Q\in {\cal
R}_{m}}D(Q||G_l), \,\,\,\,\,\,l=\overline
{1,M},\,\,\,m\not=l,\,\,\,m=\overline {1,M-1},
\end {equation}
\begin {equation}
\label {13}
E^*_{M|m}=E^*_{M|m}(E_{1|1},E_{2|2},...,E_{M-1|M-1})\df\inf\limits_{P\in
{\cal R}_{M}}D(Q||G_m), \,\,\,\,\,\, m=\overline {1,M-1},\end
{equation}
\begin {equation}
\label {14}
E^*_{M|M}=E^*_{M|M}(E_{1|1},E_{2|2},...,E_{M-1|M-1})\df\min\limits_{m:m=\overline
{1,M-1}}E^*_{M|m}.
\end {equation}

{\bf Theorem 3.1 \cite {E5}:} {\it If for described model all
conditional PDs $G_m$, $m=\overline {1,M}$,  are different in the
sense  that, $D(G_l||G_m)>0$, $l\not=m$,
  and the positive numbers
$E_{1|1},E_{2|2},...$, $E_{M-1|M-1}$ are such that the following
$M-1$ inequalities, called compatibility conditions,  hold
$$
E_{1|1}  <\min\limits_{m=\overline {2,M}}D(G_m||G_1),
$$
\vspace {-8mm}
\begin {equation}
\label {15}
\end {equation}
$$E_{m|m} <\min \left[\min\limits_{l=\overline
{1,m-1}}E^*_{l|m}(E_{l|l}),\,\,\, \min\limits_{l=\overline
{m+1,L}}D(G_l||G_m)  \right], \,\, m=\overline {2,M-1},$$ then there
exists a $LAO$ sequence $\varphi^*$ of tests, the reliabilities
matrix of which ${\bf E}(\varphi^*)=\{E^*_{m|l}\}$ is defined in
$(\ref {11})$--$(\ref {14})$ and all elements of it are positive.

When  one of   inequalities $(\ref {15})$ is violated, then  at
least one element of  matrix $\bf E(\varphi^*)$ is equal to $0$}.

The proof of Theorem 3.1 is postponed to the Appendix.

It is worth to formulate the following useful property of
reliabilities matrix of the LAO test.

 {\bf Remark 3.1 \cite {My6}:} {\it The
diagonal elements of the reliabilities matrix of the LAO test in
each row are equal only to the element of   the last column:
\begin{equation}
E_{m|m}^*=E_{M|m}^*,\,\,\, \mbox {\em and}\,\,\,E_{m|m}^*<
E_{l|m}^*,\,\,\,l =\overline {1,M-1},\,\,\,l\not=m, \,\,\, m
=\overline {1,M}.
 \label{16}
\end{equation}
That is the elements of the last column are equal to the diagonal
elements of the same row and due to $(\ref {7})$ are minimal in this
row. Consequently the first $M-1$ elements of the last column also
can play a part as  given  parameters for construction of a LAO
test}.

\section {Reliability Approach to Multihypotheses Testing for Many  Objects}

In   \cite {R}  Ahlswede and Haroutunian proposed  a new aspect of
 the statistical theory -- investigation of models with  many objects.
  This work developed  the ideas of papers on
 Information theory
 \cite {R1},  \cite {A2}, of   papers on many hypotheses testing  \cite {E4}-\cite {E5}   and  of book \cite
 {B1},
  devoted  to  research of sequential
  procedures  solving decision problems such as ranking and identification.  The problem of
hypotheses testing for the model consisting  of two independent and
of two strictly dependent objects (when they cannot admit the same
distribution) with  two possible hypothetical distributions were
solved in \cite {R}. In \cite {My6}   the specific characteristics
of the model consisting of $K(\geq 2)$ objects each independently of
others following  one of given $M (\geq 2)$ probability
distributions were explored. In \cite {Ear}  the model composed by
stochastically related objects was investigated. The result
concerning   two independent Markov chains is presented  in \cite
{Eng}. In this section we expose these results.

\subsection {Multihypotheses LAO Testing for Many Independent Objects}

Let us now consider the model with   three independent similar
objects.   For brevity we solve the problem for three objects, the
generalization of the problem  for $K$ independent objects will be
discussed hereafter along  the text.

Let $X_1$, $X_2$ and $X_3$ be independent  RVs taking values in the
same finite set ${\cal X}$, each of them  with one of $M$
hypothetical PDs $G_m=\{G_m(x),\,\,\,\,x\in {\cal X}\}$. These RVs
are the characteristics of the  objects. The random vector
$(X_1,X_2,X_3)$ assumes values $(x^1,x^2,x^3)\in {\cal X}^3$.

 Let    $({\bf x_1},{\bf x_2},{\bf x_3})\df $
$((x_{1}^{1},x_{1}^{2},x_{1}^{3}),...\, ,
(x_{n}^{1},x_{n}^{2},x_{n}^{3}),...\, ,
(x_{N}^{1},x_{N}^{2},x_{N}^{3}))$, $x_n^k\in {\cal X},$ $k=\overline
{1,3}$,
 $n=\overline {1,N},$  be a vector  of results of $N$ independent observations of the
 family   $(X_1,X_2,X_3)$.
 The test has  to determine   unknown PDs of the objects  on the base of   observed data.
  The detection for each object should  be
 made from   the same set of hypotheses: $G_m$, $m=\overline {1,M}$.
We call  this  procedure  the  {\it compound test} for three objects
and denote it by $\Phi_N$, it can be composed of   three individual
tests $ \varphi^{1}_N$, $\varphi^{2}_N$, $\varphi^{3}_N$ for each of
three  objects. The test $\varphi^{i}_N$, $i=\overline {1,3}$, is a
division of the space  ${\cal X}^N$ into $M $ disjoint subsets
${\cal A}_m^i$, $m=\overline {1, M}$. The set ${\cal A}_m^i$,
$m=\overline {1,M}$,  contains all vectors ${\bf x}_i$ for which the
hypothesis $G_m$ is adopted. Hence test  $\Phi_N$ is realised by
division of the space ${\cal X}^N\times {\cal X}^N\times {\cal X}^N$
into $ M^3$ subsets ${\cal A}_{m_1,m_2,m_3}={\cal A}_{m_{1}}^1\times
{\cal A}_{m_2}^2\times {\cal A}_{m_3}^3$, $m_i=\overline {1,M}$,
$i=\overline {1,3}$.
  We denote the infinite sequence of compound tests by
${{\Phi}}$. When we have $K$ independent objects the  test
${{\Phi}}$ is composed of  tests $ \varphi^{1}$, $\varphi^{2}$,\,
... \, , $\varphi^{K}$.

  The probability of the falsity of  acceptance    of hypotheses triple
$(G_{l_{1}},G_{l_{2}},G_{l_{3}})$ by the test $\Phi_N$ provided that
the triple of hypotheses $(G_{m_{1}},G_{m_{2}},G_{m_{3}})$ is true,
where $(m_1,m_2,m_3)\not=(l_1,l_2,l_3)$, $m_i,l_i= \overline {1,M}$,
$i=\overline {1,3}$,  is:
\begin {align*}
\alpha_{l_1,l_2,l_3|m_1,m_2,m_3}(\Phi_N)& \df G_{m_1}^N\circ
G_{m_2}^N\circ G_{m_3}^N \left ({\cal
A}^N_{l_1,l_2,l_3}\right )\nonumber\\
& \df G_{m_1}^N\left ({\cal A}^N_{l_1}\right )\cdot G_{m_2}^N\left
({\cal A}^N_{l_2}\right )\cdot
G_{m_3}^N\left ({\cal A}^N_{l_3}\right )\nonumber\\
&=\sum\limits_{{\bf x}_1\in{\cal A}^N_{l_1} }G_{m_1}^N({\bf
x}_1)\sum\limits_{{\bf x}_2\in{\cal A}^N_{l_2} }G_{m_2}^N({\bf
x}_2)\sum\limits_{{\bf x}_3\in{\cal A}^N_{l_3} }G_{m_3}^N({\bf
x}_3).\nonumber
 \end {align*}

 The probability to reject a true triple of
hypotheses $(G_{m_{1}},G_{m_{2}},G_{m_{3}})$ by analogy with (\ref
{5}) is defined as  follows:
\begin{equation}
\alpha
_{m_1,m_2,m_3|m_1,m_2,m_3}(\Phi_N)\df\sum\limits_{(l_1,l_2,l_3)\not=(m_1,m_2,m_3)}
\alpha_{l_1,l_2,l_3|m_1,m_2,m_3}(\Phi_N).
 \label{17}
\end{equation}

We  study corresponding  reliabilities
  $E_{l_1,l_2,l_3|m_1,m_2,m_3}(\Phi)$   of the sequence of
tests $\Phi$,
$$
E_{l_1,l_2,l_3|m_1,m_2,m_3}(\Phi)\df\overline{\lim\limits_{N\to\infty}}\left
\{-\frac{1}{N} \log\alpha_{l_1,l_2,l_3|m_1,m_2,m_3} (\Phi_N)\right
\},
$$
\begin{equation}
 m_i,l_i= \overline {1,M},\,\,\,i=\overline {1,3}.
 \label{18}
\end{equation}

 Definitions (\ref {17}) and (\ref {18})  imply (compare with (\ref {7})) that
\begin{equation}
E_{m_1,m_2,m_3|m_1,m_2,m_3}(\Phi)=\min\limits_{(l_1,l_2,l_3)\not=(m_1,m_2,m_3)}E_{l_1,l_2,l_3|m_1,m_2,m_3}(\Phi).
 \label{19}
\end{equation}

Our aim is  to analyze the   reliabilities matrix ${\bf
E}(\Phi^*)=\{E_{l_1,l_2,l_3|m_1,m_2,m_3}(\Phi^*)\}$ of LAO test
sequence $\Phi^*$ for three objects. We call the   test sequence LAO
for the model with many  objects if for given positive values of
 certain part of  elements  of  reliabilities  matrix
 the procedure provides maximal values for all other  elements of it.

 Let us
denote by  ${\bf E}{(\varphi^i)}$ the reliabilities  matrices of the
sequences of tests  $\varphi^i$, $i=\overline {1,3}$. The following
Lemma is a generalization of Lemma from \cite {R}.

{\bf Lemma 4.1:} {\it If    elements $E_{l|m}{(\varphi^i)}$,
$m,l=\overline {1,M}$, $i=\overline {1,3},$ are  strictly
 positive, then the following equalities hold for ${\bf E}(\Phi)$,
$\Phi=(\varphi^1,\varphi^2,\varphi^3)$,  $l_i,m_i=\overline {1,M}$:
\begin{align}
 E_{l_1,l_2,l_3|m_1,m_2,m_3}(\Phi)& =  \sum\limits_{i=1}^3E_{l_i|m_i}(\varphi^i),\nonumber\\
 \label{20} &  \,\,\, m_i\not=l_i,\\
 E_{ l_1,l_2,l_3|m_1,m_2,m_3}(\Phi)   &=  \sum\limits_{i\in
[[1,2,3]-k]}E_{l_i|m_i}{(\varphi^i)},\nonumber\\
 \label{21} &  m_{k}=l_{k},\, m_i\not=l_i, \,\,k=\overline {1,3},\\
 E_{ l_1,l_2,l_3|m_1,m_2,m_3}(\Phi) &=
 E_{l_i|m_i}{(\varphi^i)},\nonumber\\
 \label{22} & i=\overline {1,3}, \,m_{k}=l_{k},\, m_i\not=l_i, \,\,k\in
 [[1,2,3]-i].\\
 \nonumber
\end {align}

Equalities $(\ref {20})$ are valid also if
$E_{l_i|m_i}(\varphi^i)=0$ for several pairs  $(m_i,l_i)$ and
several $i$.}

 The proof of Lemma 4.1 is exposed in Appendix.

 Now we shall show how we can  erect the LAO
test from the set of compound tests when $3(M-1)$ strictly positive
elements  of the reliabilities matrix $E_{M,M,M|m,M,M}$,
$E_{M,M,M|M,m,M}$ and $E_{M,M,M|M,M,m}$, $m=\overline {1,M-1}$, are
preliminarily given.

The following subset of tests:
$$
{\cal D}=\{ \Phi: E_{m|m}(\varphi^i)> 0,\,\,\,\,m=\overline
{1,M},\,\,\,i=\overline {1,3}\}
$$
is\, distinguished\,\, by \,the\, property \,that  \, when
\,\,$\Phi\in {\cal D}$\,\, the \, elements\\
$E_{M,M,M|m,M,M}(\Phi)$, $E_{M,M,M|M,m,M}(\Phi)$ and
$E_{M,M,M|M,M,m}(\Phi)$, $m=\overline {1,M-1}$, of the reliabilities
matrix  are strictly positive.

Really, because $E_{m|m}(\varphi^i)> 0$, $m=\overline {1,M}$,
$i=\overline {1,3}$, then in view of (\ref {7}) $E_{M|m}(\varphi^i)$
are also strictly positive. From equalities (\ref {20})--(\ref {22})
we obtain that the noted elements are strictly positive for
$\Phi\in{\cal D}$ and $m=\overline {1,M-1}$
\begin{equation}
 \label{23_1} E_{M,M,M|m,M,M}(\Phi)=E_{M|m}(\varphi^1),
 \end{equation}
 \begin{equation}
\label{23_2}E_{M,M,M|M,m,M}(\Phi)=E_{M|m}(\varphi^2),
\end{equation}
\begin{equation}
\label{23_3}E_{M,M,M|M,M,m}(\Phi)=E_{M|m}(\varphi^3).
  \end{equation}

 For given   positive elements
 $$E_{M,M,M|m,M,M}, \,\,E_{M,M,M|M,m,M}, \,\,E_{M,M,M|M,M,m},  \,\, m=\overline {1,M-1},$$  define
 the  following family of decision  sets of PDs:
$$
{\cal R}_m^{(i)}\df \{Q:D(Q||G_m)\leq E_{M,M,M|m_1,m_2,m_3},
\,\,m_i=m, \,\,m_j=M,\,\,\,i\not=j,j=\overline {1,3}\}
$$
\begin{equation}
 m=\overline
{1,M-1},\,\,\, i=\overline {1,3}, \label{n1}
\end{equation}
$$
{\cal R}_{M}^{(i)}\df \{Q: D(Q||G_m)>
E_{M,M,M|m_1,m_2,m_3},\,\,m_i=m, \,\,m_j=M,\,i\not=j,j=\overline
{1,3},
$$
\begin{equation}
m=\overline {1,M-1}\}, \,\,\,i=\overline {1,3}.
 \label{n2}
\end{equation}

Define also the elements  of the reliability matrix of the compound
LAO test for three objects:
\begin{align}
E_{M,M,M|m,M,M}^*  &\df E_{M,M,M|m,M,M},
\nonumber\\
\label{24} E_{M,M,M|M,m,M}^*  &\df E_{M,M,M|M,m,M},\\
E_{M,M,M|M,M,m}^* & \df E_{M,M,M|M,M,m},\nonumber\\
 E_{ l_1,l_2,l_3|m_1,m_2,m_3}^*  &\df \inf\limits_{Q\in
R_{l_i}^{(i)}}D(Q||G_{m_i}),\nonumber\\
 \label{25}  &i=\overline {1,3}, \, m_{k}=l_{k},\,m_i\not=l_i,\,i\not=k,\,k\in [[1,2,3]-i],\\
 E_{ m_1,m_2,m_3|l_1,l_2,m_3}^*  &\df \sum\limits_{i\not=
k}\inf\limits_{Q\in R_{l_i}^{(i)}}D(Q||G_{m_i}),\nonumber\\
\label{26}  & m_{k}=l_{k},\,m_i\not=l_i,\,k=\overline {1,3},\,i\in
[[1,2,3]-k],\\
\label{27}   E_{l_1,l_2,l_3|m_1,m_2,m_3}^*  &\df \sum\limits_{i=
1}^3\inf\limits_{Q\in R_{l_i}^{(i)}}D(Q||G_{m_i}),\, m_i\not=l_i,
i=\overline {1,3}.\\
\nonumber
 \end{align}

The following theorem   is a generalization and improvement of the
corresponding theorem proved in \cite {R} for the case $K=2,
\,\,M=2$.

 {\bf Theorem 4.1 \cite {My6}:} {\it For considered model with three objects, if all
 distributions  $G_m$, $m=\overline {1,M}$,  are
different, (and equivalently $D(G_l||G_m)>0$,  $l\not=m$,
$l,m=\overline {1,M}$), then the following    statements are valid:

a) when   given strictly positive elements $E_{M,M,M|m,M,M}$,
 $E_{M,M,M|M,m,M}$ and $E_{M,M,M|M,M,m}$, $m=\overline {1,M-1}$,  meet  the following conditions
 $$
 \hspace {-2cm}\max ( E_{M,M,M|1,M,M},E_{M,M,M|M,1,M},  E_{M,M,M|M,M,
1})
 $$
\begin{align}
  \label{28}& <
\min\limits_{l= \overline { 2,M }}D(G_l||G_1),\\
&  \hspace {-1cm}\mbox {and \,\,for}\,\,\,   m=\overline {2,M-1},\nonumber\\
 \label{29} E_{M,M,M |m,M,M} &<  \min\left [\min\limits_{l= \overline {1,m-1 }}
E_{l,m,m|m,m,m}^*,\,\min\limits_{l= \overline{
m+1,M}}D(G_l||G_m)\right ], \\
\label{30} E_{M,M,M|M,m,M} &< \min\left [\min\limits_{l= \overline
{1,m-1}} E_{m,l,m|m,m,m}^*,\,\min\limits_{l= \overline { m+1,M
}}D(G_l||G_m)\right ], \\
\label{31} E_{M,M,M|M,M,m} &<  \min\left [\min\limits_{l= \overline
{1,m-1}} E_{m,m,l|m,m,m}^*,\,\,\min\limits_{l= \overline { m+1,M
}}D(G_l||G_m)\right ],\\
 \nonumber
\end{align}
then there exists a LAO test sequence $\Phi^*\in {\cal D}$, the
reliability matrix of  which
  ${{\bf E}(\Phi^*)}$ is defined in $(\ref {24})$--$(\ref
  {27})$
  and all elements of
  it are positive,

b) if even  one of the inequalities $(\ref {28})$--$(\ref {31})$ is
violated, then there exists at least one element of  the matrix
 ${{\bf E}(\Phi^*)}$ equal to $0$.}

For the proof of Theorem 4.1 see Appendix.

 When we consider the model with $K$ independent objects the
 generalization of Lemma 4.1 will take  the  following form.

 {\bf Lemma 4.2:} {\it   If  elements $E_{l_i|m_i}{(\varphi^i)}$,
  $m_i,l_i=\overline {1,M}$, $i=\overline {1,K},$ are  strictly positive,
 then the following equalities hold for}
 $\Phi=(\varphi^1,\varphi^2,..., \varphi^K)$:
\begin {align*}
E_{l_1,l_2,...,l_K
 |m_1,m_2,...,m_K}(\Phi)&=  \sum\limits_{i=1}^KE_{l_i|m_i}(\varphi^i),\,\,\,
\,\,\, m_i\not=l_i,\, i=\overline {1,K},\nonumber\\
E_{l_1,l_2,...,l_K|m_1,m_2,...,m_K}(\Phi)&=
\sum\limits_{i:\,\,m_i\not=l_i}E_{l_i|m_i}(\varphi^i).\nonumber
\end {align*}

 For  given $K(M-1)$  strictly
  positive elements $E_{M,M,...,M|m,M,...,M}$,\\ $E_{M,M,...,M|M,m,...,M}$, .... , $E_{M,...,M,M,|M,M...,m}$,
  $m=\overline {1, M-1}$, for $K$ independent objects we can find  the LAO test
  $\Phi ^*$ in a way similar to  case of three independent objects.

 {{\bf Comment 4.1: }
 Idea to  renumber
   $K$-distributions  from $1$ to $M^K$ and consider them as PDs
of one complex object offers  an alternative way  of testing for
models with  $K$ objects. We can give $M^K-1$ diagonal elements of
corresponding large  matrix ${\bf E}(\Phi)$ and apply Theorem 3.2
concerning one composite  object. In this direct algorithm  the
number of the preliminarily given elements of the matrix ${\bf E}
(\Phi)$ would be greater (because $M^K-1 > K(M-1)$, $M\geq 2, K\geq
2$) but the procedure of calculations would be longer than in our
algorithm presented above in this section. Our    approach to the
problem gives also  the possibility to define the LAO tests for each
of the separate objects, but the approach with renumbering of
$K$-vectors of hypotheses does not have this opportunity. In the
same time in the case of direct algorithm  there is opportunity for
the investigator to define preliminarily  the greater number of
elements of the matrix ${\bf E}(\Phi)$.

In applications   one of  two approaches may be used in
 conformity with preferences of the investigator.}

\subsubsection {Example}

Some illustrations of exposed results are in an  example concerning
two objects. The set ${\cal X}=\{0,1\}$ contains two elements and
the following PDs are given on ${\cal X}$:\,\,\, $G_1=\{0,10;
0,90\},$ $G_2=\{0,85; 0,15\},$ $G_3=\{0,23; 0,77\}$. As it follows
from relations (\ref {24})--(\ref {27}),  several elements of the
reliability matrix are functions of one of given elements, there are
also elements which are functions of two, or three given elements.
For example,   in Fig. 1 and Fig. 2 the results of calculations of
functions $E_{1,2|2,1}(E_{3,3|1,3},E_{3,3|3,2})$ and
$E_{1,2|2,2}(E_{3,3|1,3})$ are presented. For  these  distributions
we  have $\min(D(G_2||G_1), D(G_3||G_1))\approx 2,2$ and
$\min(E_{2,2|2,1}, D(G_3||G_2))\approx 1,4$. We see that,  when the
inequalities  (\ref {28}) or  (\ref {31}) are violated,
$E_{1,2|2,1}=0$ and $E_{1,2|2,2}=0$.
\newpage
\begin{figure}
\includegraphics[width=3.5in]{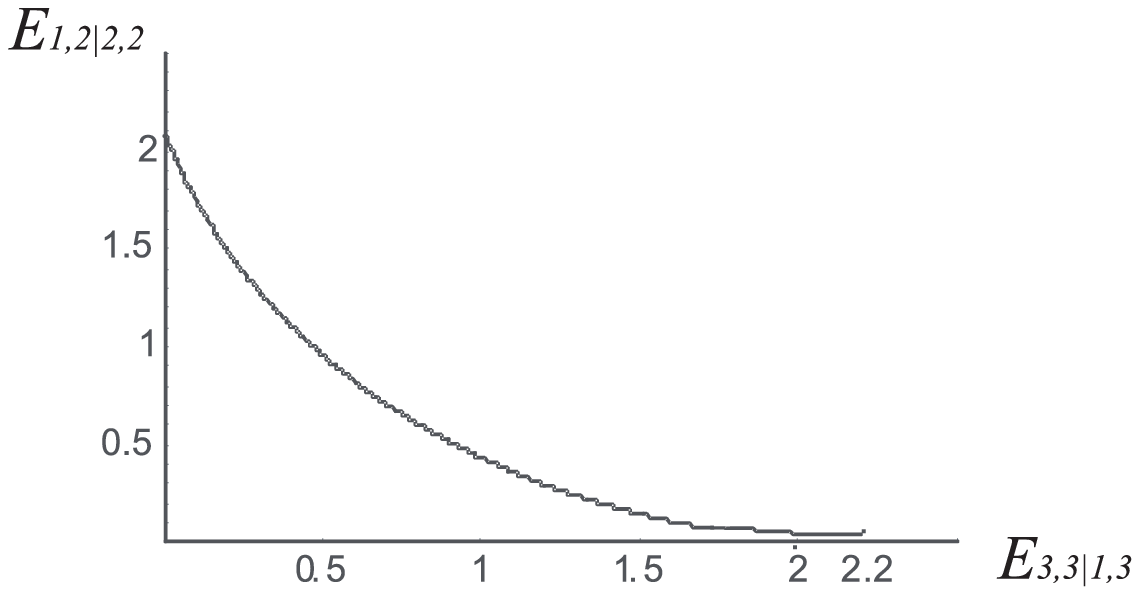}
\vspace{0.5pc} \caption []{} \label{penG}
\end{figure}

\begin{figure}
\includegraphics[width=3.5in]{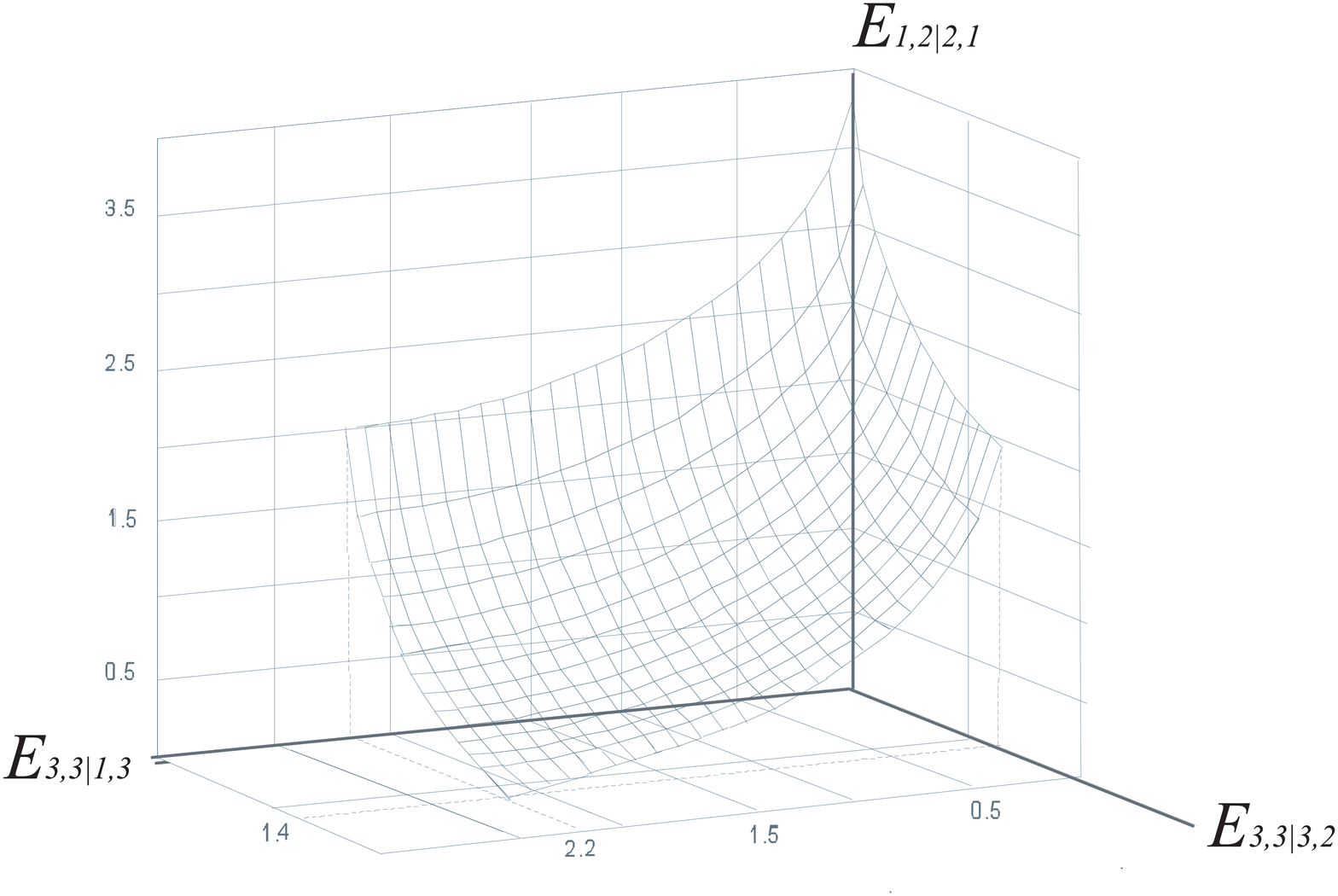}
\vspace{0.5pc} \caption []{} \label{penG1}
\end{figure}

\newpage

\subsection {Multihypotheses LAO Testing for Two  Dependent Objects}

We consider  characteristics of procedures of
 LAO testing of probability distributions  of  two related
 ({\it stochastically}, {\it statistically} and {\it strictly dependent})
  objects. We use these terms for different kinds of dependence of two
  objects.

Let $X_1$ and $X_2$ be  RVs  taking values in the same finite set
${\cal X}$ and ${\cal P}({\cal X})$ be the space of all possible
distributions on ${\cal X}$.

Let $({\bf x_1,x_2})=((x_1^1,x_1^2),(x_2^1,x_2^2),...(x_N^1,x_N^2))$
 be a sequence of results of $N$ independent observations of the pair of objects.

First we consider the model, which consists of two  {\it
stochastically} related objects. We name so the following more
general dependence.  There are given $M_1$ PDs
$$G_{m_1}=\{G_{m_1}(x^1),\,\,{x^1}\in {\cal X}\},\,\,\,
m_1=\overline{1,M_1}.$$
 The first object
is characterized by RV  $X_1$ which has one of these $M_1$ possible
PDs and the second object is dependent on the first and is
characterized by RV $X_2$ which can have one of  $M_1\times M_2$
conditional PDs
$$G_{m_2|m_1}=\{G_{m_2|m_1}(x^2|x^1),\,\,{x}^1,
{x}^2\in {\cal X}\},\,m_1=\overline{1,M_1},\,m_2=\overline{1,M_2}.
$$
Joint PD of the pair of objects is
$$G_{m_1,m_2}=G_{m_1} \circ  G_{m_2|m_1}=\{G_{m_1,m_2}({ x^1,x^2}),\,\,{x}^1,
{x}^2\in {\cal X}\},  $$ where
$$
G_{m_1,m_2}({ x^1,x^2})=G_{m_1}({ x^1})G_{m_2|m_1}({x^2}|{ x^1}),\,\,\,
m_1=\overline{1,M_1},\,m_2=\overline{1,M_2}.
$$
The probability $G_{m_1,m_2}^N({\bf x_1,x_2})$ of $N$-vector $({\bf
x_1,x_2})$  is  the following product:
\begin {align*}
G_{m_1,m_2}^N({\bf x_1,x_2})&\df G_{m_1}^N({\bf
x_1})G_{m_2|m_1}^N({\bf x_2|x_1})\nonumber\\
&\df\prod\limits_{n=1}^N
G_{m_1}(x_n^1)G_{m_2|m_1}(x_n^2|x_n^1),\nonumber
 \end {align*}
 with $G_{m_1}^N({\bf x_1})=\prod\limits_{n=1}^NG_{m_1}(x_n^1)$ and
$G_{m_2m_1}^N({\bf
x_2|x_1})=\prod\limits_{n=1}^NG_{m_2|m_1}(x_n^2|x_n^1).$

In  somewhat  particular  case, when $X_1$ and $X_2$ are related
{\it statistically} \cite {A6}, \cite {A7}, the second object
depends on index of  PD of the first object but not depends on value
$x^1$ taken by the first object. The second object is characterized
by RV $X_2$ which can have one of $M_1\times M_2$ conditional PDs
$G_{m_2|m_1}=\{G_{m_2|m_1}(x^2),\,\, {x}^2\in {\cal
X}\},\,m_1=\overline{1,M_1},\,m_2=\overline{1,M_2}$.

In the  third case  of {\it strict} dependence,  the objects $X_1$
and $X_2$ can have only different distributions from the same given
family of $M$
 PDs $G_1$, $G_2$, ..., $G_M$.

Discussed in Comment 4.1  the {\it direct  approach} for  LAO
testing of PDs for two related objects,  consisting in consideration
of the pair of objects as one composite object and then use of
Theorem 3.1, is applicable for first two cases \cite {A8}. But now
we consider also  { another approach}.

Let us remark that test $\Phi^N$ can be composed of a pair of tests
$\varphi_1^N$ and $\varphi_2^N$ for the separate objects:
$\Phi^N=(\varphi^N_1, \varphi^N_2)$. Denote by $\varphi^1$,
$\varphi^2$ and $\Phi$ the infinite sequences of tests for the
first,  for the second and for the pair of objects, respectively.

Let $X_1$ and $X_2$ be related {\it stochastically}. For the  object
characterized by $X_1$ the non-randomized test $\varphi_N^1({\bf
x_1})$ can be determined by partition of the sample space ${\cal
X}^N$   on $M_1$ disjoint subsets
 ${\cal A}^N_{l_1}=\{{\bf
x_1:}\,\, \varphi^1_N({\bf x_1})=l_1\},\,l_1=\overline{1,M_1}$,
 i.e. the  set  ${\cal
A}^N_{l_1}$ consists of  vectors $\bf x_1$ for which the PD
$G_{l_1}$ is adopted. The probability $ \alpha_{l_1|m_1}
(\varphi^1_N)$ of the erroneous acceptance of PD $G_{l_1}$  provided
that $G_{m_1}$ is true, $l_1,m_1=\overline{1,M_1},\,\,\,m_1\not=
l_1,$ is defined by the set ${\cal A}_{l_1}^{N}$
$$\alpha_{l_1|m_1} (\varphi^1_N)\df G_{m_1}^N({\cal A}_{l_1}^{N}).$$

We define  the probability  to reject $G_{m_1}$, when it is true, as
follows
\begin {equation}
\label {32} \alpha_{m_1|m_1}(\varphi^1_N)\df
\sum\limits_{l_1:l_1\not=m_1}
\alpha_{l_1|m_1}(\varphi^1_N)=G_{m_1}^N(\overline {{\cal
A}^N_{m_1}}).
\end {equation}

Corresponding error probability exponents     are:
\begin {equation}
\label {33} E_{l_1|m_1}(\varphi^1)\df \overline {\lim\limits_{N\to
\infty}}\left \{-\frac{1}{N}\log\alpha_{l_1|m_1} (\varphi^1_N)\right
\},\,\,\,m_1,l_1=\overline{1,M_1}.
\end {equation}
 It follows from (\ref {32}) and (\ref {33})
that
$$E_{m_1|m_1}(\varphi^1)=\min_{l_1:l_1\not=m_1}E_{l_1|m_1}(\varphi^1),
\,\,\,l_1,m_1=\overline{1,M_1}.$$

For construction of the LAO test we assume  given strictly positive
numbers $E_{m_1|m_1}$, $m=\overline {1,M_1-1}$ and we define regions
${\cal R}_{l_1}$, $l=\overline {1,M_1}$ as in (\ref {8})--(\ref
{9}).

For the second object characterized by RV  $X_2$ depending on $X_1$
the non-randomized test $\varphi^2_N({\bf x_2}, {\bf x_1}, l_1)$,
based on vectors $({\bf x_1, x_2})$ and on the index of the
hypothesis $l_1$ adopted for  $X_1$,  can be  given for each $l_1$
and ${\bf x}_1$ by division of the sample space ${\cal X}^N$ on
$M_2$ disjoint subsets
\begin {equation}
\label {new1}
 {\cal A}^N_{l_2|l_1}({\bf
x_1})\df \{{\bf x_2:}\,\,\, \varphi_2^N({\bf x_2}, {\bf x_1},
l_1)=l_2\},\, l_1=\overline{1,M_1},\,l_2=\overline{1,M_2}.
\end {equation}

 Let
\begin {equation}
\label {new2} {\cal A}_{l_1,l_2}^{N}\df \{ ({\bf x_1}, {\bf x_2}):
\,{\bf x_1}\in {A_{l_1}^N},\,\,{\bf x_2}\in{\cal A}^N_{l_2|l_1}({\bf
x_1})\}.
\end {equation}

The probabilities of the erroneous acceptance for $(l_1,l_2)\ne
(m_1,m_2)$ are $$ \alpha_{l_1,l_2|m_1,m_2}\df G_{m_1,m_2}^N({\cal
A}_{l_1,l_2}^{N}).
$$
Corresponding reliabilites are denoted  $E_{l_1,l_2|m_1,m_2}$ and
are  defined as in  (\ref {18}). We can  upper estimate the
probabilities of the erroneous acceptance for $(l_1,l_2)\ne
(m_1,m_2)$
\begin {align*}
  G_{m_1,m_2}^N({\cal A}_{l_1,l_2}^{N})&=\sum\limits_ {({\bf x}_1,{\bf x}_2)
  \in {\cal A}_{l_1,l_2}^N}G_{m_1}^N
 ({\bf x}_1)G_{m_2|m_1}^N({\bf x}_2|{\bf x_1})\nonumber\\
  &= \sum\limits_ {{\bf x}_1
 \in {\cal A}_{l_1}^N}G_{m_1}^N({\bf x}_1)G_{m_2|m_1}^N({\cal A}_{l_2|l_1}^{N}
 ({\bf x}_1)|{\bf x_1})\nonumber\\
 &\leq  \max\limits_{{\bf x}_1\in {\cal A}_{l_1}^N} G_{m_2|m_1}^N({\cal A}_{l_2|l_1}^{N}
 ({\bf x}_1)|{\bf x_1})\sum\limits_{{\bf x}_1\in {\cal A}_{l_1}^N}G_{m_1}^N({\bf
 x}_1)\nonumber\\
 &=  G_{m_1}^N({\cal A}_{l_1}^{N})\max\limits_{{\bf x}_1
 \in {\cal A}_{l_1}^N} G_{m_2|m_1}^N({\cal A}_{l_2|l_1}^{N}
 ({\bf x}_1)|{\bf x_1}).\nonumber
 \end {align*}

These  upper estimates of $\alpha_{l_1,l_2|m_1,m_2}(\Phi_{N})$ for
each $(l_1,l_2)\ne (m_1,m_2)$
 we denote by
\begin  {align*}
\beta_{l_1,l_2|m_1,m_2}(\Phi_{N}) &\df  G_{m_1}^N ({\cal
A}_{l_1}^{N})
 \max\limits_{{\bf x}_1\in {\cal A}_{l_1}^N}
 G_{m_2|m_1}^N({\cal A}_{l_2|l_1}^{N}({\bf x}_1)|{\bf x_1}).
  \end {align*}

  Consequently we can deduce that for $
l_1,m_1=\overline{1,M_1},\,\, l_2,m_2=\overline{1,M_2}, $ new
parameters
$$
F_{l_1,l_2|m_1,m_2}(\Phi)\df\overline{\lim\limits_{N\to\infty}}
\{-\frac{1}{N}\log\beta_{l_1,l_2|m_1,m_2}^{N} (\Phi_N)\},$$
are lower estimates for reliabilities $E_{l_1,l_2|m_1,m_2}(\Phi).$\\
We can  introduce for
$l_1,m_1=\overline{1,M_1},\,l_2,m_2=\overline{1,M_2},\,\,m_2\not=
l_2,$
$$\beta_{l_2|l_1,m_1,m_2}(\varphi^2_{N})\df\max\limits_{{\bf
x}_1\in {\cal A}_{l_1}^N} G_{m_2|m_1}^N ({\cal A}_{l_2|l_1}^N({\bf
x}_1)|{\bf x_1}),$$
 and  also consider
$$\beta_{m_2|l_1,m_1,m_2}(\varphi^2_{N})\df\max\limits_{{\bf x}_1\in {\cal A}_{l_1}^N}
G_{m_2|m_1}^N(\overline{{\cal A}_{m_2|l_1}^N({\bf x}_1)}|{\bf x_1})
$$
\begin {equation}
\label {34}\hspace{14mm}=
\sum\limits_{l_2\not=m_2}\beta_{l_2|l_1,m_1,m_2}(\varphi^2_{N}).
\end {equation}
 The corresponding  estimates of the reliabilities of test
$\varphi^2_{N}$, are the following
$$
F_{l_2|l_1,m_1,m_2}(\varphi^2)\df \overline {\lim\limits_{N\to
\infty}}\left \{-\frac{1}{N}\log\beta_{l_2|l_1,m_1,m_2}
(\varphi^2_{N})\right \},
$$
\begin {equation}
\label {35} l_1,m_1=\overline{1,M_1},\, l_2,m_2=\overline{1,M_2},\,
m_2\not= l_2.
\end {equation}
 It is
clear from (\ref {34}) and (\ref {35}) that for every
$l_1,m_1=\overline{1,M_1}, \, l_2,m_2=\overline{1,M_2}$
\begin {equation}
\label {36}
F_{m_2|l_1,m_1,m_2}(\varphi^2)=\min_{l_2:l_2\not=m_2}F_{l_2|l_1,m_1,m_2}(\varphi^2).
\end {equation}
For given  positive numbers
$F_{l_2|l_1,m_1,l_2},\,l_2=\overline{1,M_2-1}$,\, for $Q\in {\cal
R}_{l_1}$ and for each pair\,\,  $l_1, m_1=\overline{1,M_1}$   let
us define the following regions and values:
\begin {align}
\label {37}{\cal R}_{l_2|l_1}(Q)&\df
\{V:\,\,\,D(V||G_{l_2|l_1}|Q)\leq
F_{l_2|l_1,m_1,l_2}\},\,\,l_2=\overline{1,M_2-1},\\
\label {38} {\cal R}_{M_2|l_1}(Q)&\df \left
\{V:\,\,\,D(V||G_{l_2|l_1}|Q)>
F_{l_2|l_1,m_1,l_2},\,\,l_2=\overline{1,M_2-1}\right \},\\
\nonumber
\end {align}
\begin {align}
\label{39} F_{l_2|l_1,m_1,l_2}^*(F_{l_2|l_1,m_1,l_2})&\df
F_{l_2|l_1,m_1,l_2},\,\,l_2=\overline {1,M_2-1},\\
 F_{l_2|l_1,m_1,m_2}^*(F_{l_2|l_1,m_1,l_2})&\df \inf\limits_{Q\in
{\cal R}_{l_1}} \inf\limits_{V\in {\cal R}_{l_2|l_1}(Q)}D(V||G_{m_2|m_1}|Q),\nonumber\\
\label {40}   &m_2=\overline {1,M_2},\,,m_2\not=l_2,\,l_2=\overline {1,M_2-1},\\
 &  \hspace {-5cm} F_{M_2|l_1,m_1,m_2}^*(F_{1|l_1,m_1,1},...,F_{M_2-1|l_1,m_1,M_2-1})\nonumber\\
 &\df  \inf\limits_{Q\in {\cal R}_{l_1}}\inf\limits_{V\in {\cal R}_{M_2|l_1}(Q)}D(V||G_{m_2|m_1}|Q),\nonumber\\
 & m_2=\overline {1,M_2-1},\\
 &\hspace {-5cm}F_{M_2|l_1,m_1,M_2}^*(F_{1|l_1,m_1,1},...,F_{M_2-1|l_1,m_1,M_2-1})\nonumber\\
\label {41*} &\df \min\limits_{l_2=\overline
{1,M_2-1}}F_{l_2|l_1,m_1,M_2}^*.\\
\nonumber
\end {align}
We denote by ${\bf F}(\varphi_2)$ the  matrix of lower estimates for
elements of matrix  ${\bf E}(\varphi_2).$

 {\bf Theorem 4.2 \cite {Ear}:} {\it If for given $m_1, l_1=\overline{1,M_1},$  all  conditional PDs
$G_{l_2|l_1}$, $ l_2=\overline {1,M_2}$,  are different in the sense
that  $D(G_{l_2|l_1}||G_{m_2|m_1}|Q)>0$, for all $Q\in {\cal
R}_{l_1}$,\, $l_2\not=m_2,\,\,m_2=\overline {1,M_2}$, when the
strictly positive  numbers $F_{1|l_1,m_1,1}$,
$F_{2|l_1,m_1,2}$,...,$F_{M_2-1|l_1,m_1,M_2-1}$ are such that the
following compatibility conditions  hold
\begin {align}
\label {43*} F_{1|l_1,m_1,1}&< \min\limits_{l_2=\overline
{2,M_2}}\inf\limits_{Q\in {\cal
R}_{l_1}}D(G_{l_2|l_1}||G_{1|m_1}|Q),\\
F_{m_2|l_1,m_1,m_2} &< \min\left (\min\limits_{l_2=\overline
{m_2+1,M_2}}\inf\limits_{Q\in {\cal
R}_{l_1}}D(G_{l_2|l_1}||G_{m_2|m_1}|Q), \right .\nonumber\\
\label {43}&  \left . \min\limits_{l_2=\overline
{1,m_2-1}}F_{l_2|l_1,m_1,m_2}^*(F_{l_2|l_1,m_1,l_2})\right ),\,\,
m_2=\overline {1,M_2-1},
 \end {align}
 then there exists a  sequence  of tests
$\varphi^{2,*}$, such that  the lower estimates are defined in
(4.29)-- (4.32) and are  strictly positive.

Inequalities $(\ref {43*})$, $(\ref {43})$ are necessary for
existence of test sequence with  matrix ${\bf F}(\varphi^{2,*})$ of
positive lower estimates having  given elements
$F_{l_2|l_1,m_1,l_2},\, l_2=\overline{1,M_2-1}$ in diagonal. }

\vspace {3mm}

Let us  define    the    following subsets  of ${\cal P}({\cal X})$
for given  strictly positive  elements  $E_{M_1,l_2|l_1,l_2}$,
$F_{l_1,M_2|l_1,l_2}$, $l_1=\overline {1,M_1-1},\,l_2=\overline
{1,M_2-1}$:
\begin {align}
{\cal R}_{l_1}&\df  \{Q:D(Q||G_{l_1})\leq E_{M_1,l_2|l_1,l_2}\},\,
l_1=\overline {1,M_1-1},\, l_2=\overline
{1,M_2-1},\\
{\cal R}_{l_2|l_1}(Q)&\df  \{V:D(V||G_{l_2|l_1}|Q)\leq
F_{l_1,M_2|l_1,l_2}\},\, l_1=\overline {1,M_1-1},\, l_2=\overline
{1,M_2-1},\\
{\cal R}_{M_1}&\df  \{Q: D(Q||G_{l_1})>
E_{M_1,l_2|l_1,l_2},\,\,l_1=\overline {1,M_1-1},\, l_2=\overline
{1,M_2-1}\},\\
{\cal R}_{M_2|l_1}(Q)&\df  \{V: D(V||G_{l_2|l_1}|Q)>
F_{l_1,M_2|l_1,l_2},\,\,l_1=\overline {1,M_1-1}, \,l_2=\overline
{1,M_2-1}\}.
\end {align}

Assume also that
\begin {align}
 \label {31.a}F_{l_1,M_2|l_1,l_2}^{*}&\df  F_{l_1,M_2|l_1,l_2},\\
E_{M_1,l_2|l_1,l_2}^{*}&\df E_{M_1,l_2|l_1,l_2},\,\,l_1=\overline
{1,M_1-1},\, l_2=\overline
{1,M_2-1},\\
\label {31.b} E_{ l_1,l_2|m_1,l_2}^{*}&\df  \inf\limits_{Q:Q\in
R_{l_1}}D(Q||G_{m_1}), \,\,\,\, m_1\not=l_1,\\
\label {31.c} F_{l_1,l_2|l_1,m_2}^{*}&\df \inf\limits_{Q\in {\cal
R}_{l_1}}\inf\limits_{V:V\in
R_{l_2/l_1}(Q)}D(V||G_{m_2/m_1}|Q),\,\,m_2\not=l_2,\\
\label {31.d} F_{ l_1,l_2|m_1,m_2}^{*}&\df
F_{m_1,l_2|m_1,m_2}^{*}+E_{l_1,m_2|m_1,m_2}^{*},
\,\,\,\,m_i\not=l_i, \,\,i=1,2,\\
\label {31.e} F_{ m_1,m_2|m_1,m_2}^{*}&\df
\min\limits_{(l_1,l_2)\not= (m_1,m_2)}F_{ l_1,l_2|m_1,m_2}^{*}.
\end {align}

{\bf Theorem 4.3 \cite {A6}:} {\it If all  PDs  $G_{m_1}$,
$m_1=\overline {1,M_1}$,  are different,  that is\\
$D(G_{l_1}||G_{m_1})>0$, $l_1\not=m_1$, $l_1,m_1=\overline {1,M_1}$,
and all  conditional PDs $G_{l_2|l_1}$,\\ $l_2=\overline {1,M_2}$,
are also different for all $l_1=\overline{1,M_1},$ in the sense
that\\ $D(G_{l_2|l_1}||G_{m_2|m_1}|Q)>0$,  $l_2\not=m_2$, then the
following statements are valid.

When   given strictly positive elements $E_{M_1,l_2|m_1,l_2}$ and
$F_{l_1,M_2|l_1,m_2}$, $m_1=\overline {1,M_1-1},\,m_2=\overline
{1,M_2-1}$,  meet  the following conditions
\begin {align}
\label {32.a}
 E_{M_1,l_2|1,l_2}&< \min\limits_{l_1=
\overline { 2,M_1}}D(G_{l_1}||G_1),\\
\label {32.b}
  F_{l_1,M_2|l_1,1}&< \min\limits_{l_2= \overline
{2,M_2}}\inf\limits_{Q\in {\cal
R}_{l_1}}D(G_{l_2|l_1}||G_{1|m_1}|Q),\\
 E_{M_1,l_2|m_1,l_2}&< \min[\min\limits_{l_1= \overline
{1,m_1-1 }} E_{l_1,l_2|m_1,l_2}^{*},\,\,\min\limits_{l_1= \overline{
m_1+1,M_1}}D(G_{l_1}||G_{m_1})],\nonumber\\
\label {32.c} &m_1=\overline {2,M_1-1},\\
  F_{l_1,M_2|l_1,m_2}&< \min[{\min\limits_{l_2= \overline
{1,m_2-1}} F_{l_1,l_2|l_1,m_2}^*,\,\,\min\limits_{l_2= \overline {
m_2+1,M_2 }}\inf\limits_{Q\in {\cal
R}_{l_1}}D(G_{l_2|l_1}||G_{m_2|m_1}|Q)}],\nonumber\\
 \label {32.d}
& m_2=\overline {2,M_2-1},\\
\nonumber
\end {align}
 then
 there  exists  a  $LAO$  test  sequence $\Phi^*$,  the matrix of
 lower estimates    of which ${{\bf F}(\Phi^*)}=\{F_{ l_1,l_2|m_1,m_2}(\Phi^{*})\}$
 is defined in (\ref {31.a})-(\ref{31.e}) and all elements of it are positive.

When even  one of the inequalities (\ref {32.a})-(\ref {32.d}) is
violated, then at least one element of  the lower estimate matrix
${{\bf F}(\Phi^*)}$ is equal to $0$}.

\vspace {2mm}
 When $X_1$ and $X_2$ are related statistically \cite
{A6}, \cite {A7}  we will have instead of (\ref {34}), (\ref {35})
${\cal A}^N_{l_2|l_1}=\{{\bf x_2:}\,\,\, \varphi_2^N({\bf x_2},
l_1)=l_2\},\, l_1=\overline{1,M_1},\,l_2=\overline{1,M_2}$, and
${\cal A}^N_{l_1,l_2}\df \{({\bf x}_1, {\bf x}_2): \,\,{\bf x}_1\in
{\cal A}_{l_1}^N,\,\,{\bf x}_2\in {\cal A}_{l_2|l_1}^N ({\bf
x}_1)\}$. In that case we have error probabilities
\begin {align}
G_{m_1,m_2}^N({\cal A}_{l_1,l_2}^{N})&\df   \sum\limits_ {({\bf
x}_1,{\bf x}_2)\in {\cal A}_{l_1,l_2}^N}G_{m_1}^N({\bf
x}_1)G_{m_2|m_1}^N({\bf x}_2)\nonumber\\
&= \sum\limits_{{\bf x}_1\in {\cal A}_{l_1}^N}G_{m_1}^N({\bf x}_1)
\sum\limits_{{\bf x}_2\in {\cal A}_{l_2|l_1}^N}G_{m_2|m_1}^N({\bf
x}_2)\nonumber\\
&= G_{m_2|m_1}^N({\cal A}^N_{l_2|l_1})G_{m_1}({\cal A}^N_{l_1}),\,\,
(l_1,l_2)\ne (m_1,m_2).\nonumber
\end {align}
For the second object  the conditional  probabilities of the
erroneous acceptance of PD $G_{l_2|l_1}$
 provided that $G_{m_2|m_1}$ is true, for   $l_1,m_1=\overline{1,M_1}$, $l_2,m_2=\overline{1,M_2}$, are
 the following
$$\alpha_{l_2|l_1,m_1,m_2}^{N}(\varphi^{2}_{N})\df  G_{m_2|m_1}^N({\cal A}^N_{l_2|l_1}),\, l_2\ne m_2. $$
The probability  to reject $G_{m_2|m_1}$, when it is true is denoted
as follows
$$
\alpha_{m_2|l_1,m_1,m_2}^{N}(\varphi^2_{N})\df
G_{m_2|m_1}^N(\overline{{\cal A}_{m_2|l_1}}) =
\sum\limits_{l_2\not=m_2}\alpha_{l_2|l_1,m_1,m_2}^{N}(\varphi^2_{N}).
$$
Thus in the conditions and in the results of Theorems 4.2 and
Theorems 4.3 instead of conditional divergences $\inf\limits_{Q\in
{\cal R}_{l_1}}D(G_{l_2|l_1}||G_{m_2|m_1}|Q)$,\\ $\inf\limits_{Q\in
{\cal R}_{l_1}}D(V||G_{m_2|m_1}|Q)$ we will have just divergences
$D(G_{l_2|l_1}||G_{m_2|m_1})$,\\ $D(V||G_{m_2|m_1})$ and in place of
$ F_{l_2|l_1,m_1,m_2}(\Phi)$, $F_{l_1,l_2|m_1,m_2}(\Phi)$,
$l_1,m_1=\overline{1,M_1}$,  $l_2, m_2=\overline{1,M_2}$, will be
$E_{l_2|l_1,m_1,m_2}(\Phi)$,\,\, $E_{l_1,l_2|m_1,m_2}(\Phi)$,
\,\,$l_1,m_1=\overline{1,M_1}$, $l_2,m_2=\overline{1,M_2}$ . And in
that case regions defined in $(\ref {37})$, $(\ref {38})$ will be
changed as follows:
\begin {align} {\cal R}_{l_2|l_1}&\df
\{V:\,\,\,D(V||G_{l_2|l_1})\leq
E_{l_2|l_1,m_1,l_2}\},\,\,\,\,l_2=\overline{1,M_2-1},\nonumber\\
 {\cal
R}_{M_2|l_1}&\df  \{V:\,\,\,D(V||G_{l_2|l_1})>
E_{l_2|l_1,m_1,l_2},\,\,\,\,l_2=\overline{1,M_2-1}\},\nonumber
\end {align}
In case of two statistically dependent objects the corresponding
regions will be
\begin {align}
{\cal R}_{l_1}&\df  \{Q:D(Q||G_{l_1})\leq
E_{M_1,l_2|l_1,l_2}\},\,\,\,l_1=\overline {1,M_1-1},\, l_2=\overline
{1,M_2-1},\nonumber\\
{\cal R}_{l_2|l_1}&\df  \{V:D(V||G_{l_2|l_1})\leq
E_{l_1,M_2|l_1,l_2}\},\,\,\,l_1=\overline {1,M_1-1},\, l_2=\overline
{1,M_2-1},\nonumber\\
{\cal R}_{M_1}&\df \{Q: D(Q||G_{l_1})>
E_{M_1,l_2|l_1,l_2},\,\,l_1=\overline {1,M_1-1},\, l_2=\overline
{1,M_2-1}\},\nonumber\\
{\cal R}_{M_2|l_1}&\df  \{V: D(V||G_{l_2|l_1})>
E_{l_1,M_2|l_1,l_2},\,\,l_1=\overline {1,M_1-1}, \,l_2=\overline
{1,M_2-1}\}.\nonumber
\end {align}
 So in this case the matrix  of reliabilities ${\bf E}(\Phi*)=\{E*_{l_1,l_2|m_1,m_2}, \,\,l_1,m_1=\overline {1,M_1},
  \,\,l_1,m_1=\overline {1,M_1}\}$, will have the following elements:
\begin {align}
  E_{l_1,M_2|l_1,l_2}^{*}&\df   E_{l_1,M_2|l_1,l_2},\nonumber\\
E_{M_1,l_2|l_1,l_2}^{*}&\df  E_{M_1,l_2|l_1,l_2},\nonumber\\
  & l_1=\overline {1,M_1-1},\, l_2=\overline {1,L_2-1},\nonumber\\
  E_{ l_1,l_2|m_1,l_2}^{*}&\df  \inf\limits_{Q:Q\in
R_{l_1}}D(Q||G_{m_1}), \,\,\,m_1\not=l_1,\nonumber\\
   E_{l_1,l_2|l_1,m_2}^{*}&\df \inf\limits_{V:V\in
R_{l_2|l_1}}D(V||G_{m_2|m_1}),\,\,\,\,\,m_2\not=l_2,\nonumber\\
   E_{ l_1,l_2|m_1,m_2}^{*}&\df
E_{m_1,l_2|m_1,m_2}^{*}+E_{l_1,m_2|m_1,m_2}^{*}, \,\, m_i\not=l_i,
\,\,i=1,2,\nonumber\\
 E_{ m_1,m_2|m_1,m_2}^{*}&\df \min\limits_{(l_1,l_2)\not=
(m_1,m_2)}E_{ l_1,l_2|m_1,m_2}^{*}.\nonumber
\end {align}

{\bf Theorem 4.4:} \cite {A6} {\it If all PDs  $G_{m_1}$,
$m_1=\overline {1,M_1}$,  are different, that is\\
$D(G_{l_1}||G_{m_1})>0$, $l_1\not=m_1$, $l_1,m_1=\overline {1,M_1}$,
and all  conditional PDs $G_{l_2|l_1}$, $l_2=\overline {1,M_2}$,
are also different for all $l_1=\overline{1,M_1},$ in the sense
that $D(G_{l_2|l_1}||G_{m_2|m_1})>0$,  $l_2\not=m_2$, then the
following statements are valid.

When   given strictly positive elements $E_{M_1,l_2|l_1,l_2}$ and
$E_{l_1,M_2|l_1,l_2}$, $l_1=\overline {1,M_1-1}$, $l_2=\overline
{1,M_2-1}$,  meet  the following compatibility conditions
\begin {align}
  E_{M_1,l_2|1,l_2}&<  \min\limits_{l_1= \overline {
2,M_1}}D(G_{l_1}||G_1),\nonumber\\
    E_{l_1,M_2|l_1,1}&< \min\limits_{l_2= \overline
{2,M_2}} D(G_{l_2|l_1}||G_{1|m_1}),\nonumber\\
  E_{M_1,l_2|m_1,l_2}&<  \min\left [\min\limits_{l_1=
\overline {1,m_1-1 }} E_{l_1,l_2|m_1,l_2}^{*},\,\,\min\limits_{l_1=
\overline{m_1+1,M_1}}D(G_{l_1}||G_{m_1})\right ],\nonumber\\
& m_1=\overline {2,M_1-1}, \nonumber\\
  E_{l_1,M_2|l_1,m_2}&<  \min\left [{\min\limits_{l_2=
\overline {1,m_2-1}} E_{l_1,l_2|l_1,m_2}^*,\,\,\min\limits_{l_2=
\overline { m_2+1,M_2 }} D(G_{l_2|l_1}||G_{m_2|m_1})}\right ],\nonumber\\
  &  m_2=\overline {2,M_2-1},\nonumber
\end {align}
 then
 there  exists  a  $LAO$ test   sequence  $\Phi^*$, the
matrix  of  which ${{\bf E}(\Phi^*)}$
 is stated above  and all elements of it are positive.

When even  one of the compatibility conditions is violated, then at
least one element of  the  matrix ${{\bf E}(\Phi^*)}$ is equal to
$0$}.

\section {Identification of Distribution  for One and for Many Objects}
In \cite{B1} Bechhofer, Kiefer, and Sobel presented  investigations
on sequential multiple-decision procedures.  This book  concerns
principally with a particular class of problems referred to as
ranking problems. Chapter 10 of the book by Ahlswede and Wegener
\cite{A1} is devoted to statistical identification and ranking.
  Problems of distribution identification and distributions
ranking for one object applying the concept of optimality developed
in \cite {9}, \cite {4},  \cite {E4}--\cite {Em3}  were solved in
\cite {R}. In papers  \cite {My7}, \cite {EAP1}, \cite {Ear} and
\cite {lider} identification problems for models composed with   two
independent, or strictly dependent  objects were investigated.

In \cite {R}, \cite {My7},  \cite {Ear} and \cite {lider} models
considered in \cite{B1} and \cite{A1}  and variations of these
models inspired by the pioneering paper  by Ahlswede and Dueck
\cite{A2}, applying the concept of optimality developed in \cite
{9}, \cite {E4}--\cite {Em3}, \cite {4}, were studied.

First we formulate the concept of  the identification  for  one
object, which was considered in \cite {R}. There are known $M\geq 2$
possible PDs, related with the object in consideration.
Identification
  gives the answer to  the question whether $r$-th PD occured,
 or not. This answer can be given on the base of a sample $\bf x$
  and by a test $\varphi ^*_N({\bf x})$.
More precisely,  identification can be considered as an answer to
the question: is result $l$ of testing algorithm   equal to $r$
(that is $l =r$), or not equal $l$ (that is $l\not = r$).

  There are two
types of error probabilities of identification for each $r=
\overline {1, M }$:
 the probability ${\alpha}_{l\neq r|m=r}{(\varphi _N)}$ to accept $l$
 different from $r$, when $r$ is in reality, and the probability
 ${\alpha}_{l =r|m\neq r}{(\varphi _N)}$ that $r$ is accepted by test $\varphi _N$, when $r$ is not correct.

The probability ${\alpha}_{l\neq r|m=r}{(\varphi _N)}$  coincides
with the error probability of testing ${\alpha}_{r|r}{(\varphi _N)}$
(see (6)) which is equal to $\sum\limits_{l: l\neq r}{\alpha}_{l|r}
{(\varphi _N)}$. The corresponding reliability $E_{l\neq
r|m=r}(\varphi )$ is equal to $E_{r|r}(\varphi)$ which satisfies the
equality (\ref {7}).

 And what is the reliability approach to identification? It is necessary  to determine
the   dependence of optimal reliability $E_{l=r|m\neq r}^*$ upon
given $E_{l\neq r|m=r}^*=E_{r|r}^*$,
 which can be assigned a value  satisfying conditions analogical to (\ref {15}).

The result from paper \cite {R} is:

{\bf  Theorem 5.1: } {\it In the case of distinct hypothetical  PDs
$G_1,G_2,...,G_M$, for a given sample $\bf x$ we define its type
$Q$, and when $Q\in {\cal R}_l^ {(N)}$ (see (\ref {8})--(\ref{10}))
we accept the hypothesis $l$. Under condition that the  a priori
probabilities of all $M$ hypotheses are positive the reliability of
such identification $E_{l=r|m\neq r}$
 for given $E_{l\neq r|m =r}=E_{r|r}$ is the following:}
$$
E_{l=r|m\neq r}(E_{r|r})
=\min\limits_{m: m\neq r}\inf\limits_{Q: D(Q\Vert G_r)\leq
E_{r|r}}D(Q\Vert G_m), \,\,\,\, r =\overline {1, M}.
$$

We can accept the supposition  of positivity of a priory
probabilities of all hypotheses with loss of generality, because the
PD which is known to have probability $0$, that is being impossible,
must not be included in the studied family.

 Now let us consider the model consisting of two independent
objects. Let hypothetical characteristics of objects $X_1$ and $X_2$
be independent   RVs  taking values in the same finite set  ${\cal
X}$ with one of $M$  PDs. Identification means that the statistician
has to answer the question whether the pair of distributions
$(r_1,r_2)$ occurred or not.  Now the procedure of testing for two
objects can be used. Let us study two types of error probabilities
for each pair $(r_1,r_2)$, $r_1,r_2=\overline {1, M}$. We denote by
$\alpha_{(l_1,l_2)\not= (r_1,r_2)|(m_1,m_2)=(r_1,r_2)}^{(N)}$ the
probability,  that pair $(r_1,r_2)$ is true, but it is rejected.
Note that this probability is equal to
$\alpha_{r_1,r_2|r_1,r_2}(\Phi_N)$.  Let $\alpha_{(l_1,l_2)=
(r_1,r_2)|(m_1,m_2)\not=(r_1,r_2)}^{(N)}$ be the probability that
$(r_1,r_2)$ is identified, when it is not correct. The corresponding
reliabilities are $E_{(l_1,l_2)\not=
(r_1,r_2)|(m_1,m_2)=(r_1,r_2)}=E_{r_1,r_2|r_1,r_2}$ and
$E_{(l_1,l_2)= (r_1,r_2)|(m_1,m_2)\not=(r_1,r_2)}$. Our aim is to
determine the dependence of   $E_{(l_1,l_2)=
(r_1,r_2)|(m_1,m_2)\not=(r_1,r_2)}$ on given
$E_{r_1,r_2|r_1,r_2}(\Phi_N)$.

Let us define for each $r$, $r=\overline {1,M}$, the following
expression:
$$
A(r)=\min\left [\min\limits_{l= \overline {1, r-1 }}D(G_l||G_{r}),
 \min\limits_{l= \overline {r+1,M}}D(G_l||G_{r})\right ].
$$

{\bf Theorem 5.2 \cite {My7}:}{\it For the  model consisting of two
independent objects  if the distributions $G_m$, $m=\overline
{1,M}$, are different and the given strictly  positive number
$E_{r_1,r_2|r_1,r_2}$ satisfy condition
$$
E_{r_1,r_2|r_1,r_2}<\min\left[A(r_1), \,\,A(r_2)\right ],
$$
 then the reliability $E_{(l_1,l_2)=(r_1,r_2)|(m_1,m_2)\not=(r_1,r_2)}$ is defined as
 follows:
$$
 E_{(l_1,l_2)= (r_1,r_2)|(m_1,m_2)\not=(r_1,r_2)}\left
(E_{r_1,r_2|r_1,r_2}\right ) $$ $$=\min\limits_{m_1\not=r_1,m_2\not=
r_2}\left
[E_{m_1|r_1}(E_{r_1,r_2|r_1,r_2}),E_{m_2|r_2}(E_{r_1,r_2|r_1,r_2})\right],
$$
 where $E_{m_1|r_1}(E_{r_1,r_2|r_1,r_2})$ and
$E_{m_2|r_2}(E_{r_1,r_2|r_1,r_2})$ are determined by $(\ref {12})$}.
\vspace {2mm}

 Now  we will present  the lower estimates of the
reliabilities for LAO identification for the {\it dependent} object
which can be then  applied for deducing the lower estimates of the
reliabilities for LAO identification of two {\it related} objects.
There exist two error probabilities for each $r_2=\overline {1,
M_2}$: the probability ${\alpha}_{l_2\neq
r_2|l_1,m_1,m_2=r_2}{(\varphi^2_N)}$ to accept $l_2$ different from
$r_2$, when $r_2$ is in reality, and the probability\\
${\alpha}_{l_2 =r_2|l_1,m_1,m_2\neq r_2}{(\varphi^2 _N)}$ that $r_2$
is accepted, when it is not correct.

The upper estimate ${\beta}_{l_2\neq r_2|l_1,m_1,m_2=r_2}{(\varphi
^2_N)}$ for  ${\alpha}_{l_2\neq r_2|l_1,m_1,m_2=r_2}{(\varphi^2_N)}$
is already known, it coincides with
${\beta}_{r_2|l_1,m_1,r_2}{(\varphi^2_N)}$ which is equal to\\
$\sum\limits_{l_2: l_2\neq r_2}{\beta}_{l_2|l_1,m_1,r_2} {(\varphi
^2_N)}$. The corresponding  $F_{l_2\neq
r_2|l_1,m_1,m_2=r_2}(\varphi^2 )$ is equal to\\
$F_{r_2|l_1,m_1,r_2}(\varphi^2)$, which satisfies the equality
(\ref{19}).

We determine the  optimal dependence of $F_{l_2=r_2|l_1,m_1,m_2\neq
r_2}^*$ upon given\\ $F_{l_2\neq r_2|l_1,m_1,m_2=r_2}^*$.

{\bf  Theorem 5.3 \cite {Ear}:} {\it In  case of distinct PDs
$G_{1|l_1},G_{2|l_1},...,G_{M_2|l_1}$, under condition that a priori
probabilities of all $M_2$ hypotheses are strictly positive, for
each $ r_2=\overline {1, M_2}$ the estimate of
$F_{l_2=r_2|l_1,m_1,m_2\neq r_2}$  for given $F_{l_2\neq
r_2|l_1,m_1,m_2 =r_2}=F_{r_2|l_1,m_1,r_2}$ is the following:}
$$
F_{l_2=r_2|l_1,m_1,m_2\neq r_2}(F_{r_2|l_1,m_1,r_2}) =
$$
$$
\min\limits_{m_2: m_2\neq r_2}\inf\limits_{Q\in {\cal R}_{l_1}}
\inf\limits_{V: D(V\Vert G_{r_2|l_1}|Q)\leq
F_{r_2|l_1,m_1,r_2}}D(V\Vert G_{m_2|m_1}|Q).$$

The result of the reliability approach to  the problem of
identification of the probability distributions for two {\it
related} objects is the following.

{\bf Theorem 5.4:} {\it If the distributions $G_{m_1}$ and
$G_{m_2|m_1}$, $m_1=\overline {1,M_1}$, $m_2=\overline {1,M_2}$, are
different and the given strictly positive number
$F_{r_1,r_2|r_1,r_2}$ satisfies the condition
\begin {align*}
E_{r_1|r_1}&< \min\left [\min\limits_{l= \overline {1, r_1-1
}}D(G_{r_1}||G_{l_1}), \min\limits_{l_1= \overline
{r_1+1,M_1}}D(G_{l_1}||G_{r_1}) \right ],
\end {align*}
 \mbox {or}\\
 \begin {align*}
 F_{r_2|l_1,m_1,r_2}&< \min \left[\inf\limits_{Q\in {\cal
R}_{l_1}}\min\limits_{l_2= \overline {1, r_2-1
}}D(G_{r_2|m_1}||G_{l_2|l_1}|Q),\right .\\
&  \left .\inf\limits_{Q\in {\cal R}_{l_1}}\min\limits_{l_2=
\overline {r_2+1,M_2}}D(G_{l_2|l_1}||G_{r_2|m_1}|Q) \right],
\end {align*}
 then the  lower estimate  $F_{(l_1,l_2)=(r_1,r_2)|(m_1,m_2)\not=(r_1,r_2)}$
 of the reliability\\
 $E_{(l_1,l_2)=(r_1,r_2)|(m_1,m_2)\not=(r_1,r_2)}$ can be calculated
 as follows
$$
 \hspace {-5cm}F_{(l_1,l_2)= (r_1,r_2)|(m_1,m_2)\not=(r_1,r_2)}\left
(F_{r_1,r_2|r_1,r_2}\right )
$$
$$
=\min\limits_{m_1\not=r_1,m_2\not= r_2}\left
[E_{r_1|m_1}(F_{r_1,r_2|r_1,r_2}),F_{r_2|l_1,m_1,m_2}(F_{r_1,r_2|r_1,r_2})\right],
$$
where $E_{r_1|m_1}(F_{r_1,r_2|r_1,r_2})$ and
$F_{r_2|l_1,m_1,m_2}(F_{r_1,r_2|r_1,r_2})$ are determined
respectively by  (\ref {12}) and (\ref {40}).}

The particular case, when $X_1$ and $X_2$ are related {\it
statistically}, was studied in \cite {A6}, \cite {A7}.

\section {Multihypotheses Testing With  Possibility of  Rejection of Decision}
This section is devoted to description of characteristics of LAO
hypotheses testing with permission of decision rejection  for the
model consisting of one or more objects.  The multiple hypotheses
testing problem with possibility of rejection of decision for
arbitrarily varying object with side information and for the model
of two or  more independent objects was examined by Haroutunian,
 Hakobyan and Yessayan \cite {EHA}, \cite {HPA11}.
  These   works ware  induced by  the paper of  Nikulin \cite {Nk} concerning
      two hypotheses testing with refusal to take decision.
 An asymptotically optimal classification, in particular
hypotheses
 testing problem with rejection of decision ware   considered by  Gutman \cite {11}.

\subsection { Many Hypothesis Testing With Rejection of Decision by
Informed Statistician for Arbitrarily Varying Object}

In this  section we consider multiple statistical hypotheses testing
with possibility of rejecting to  make choice between hypotheses
concerning distribution of a discrete arbitrarily varying object.
The arbitrarily varying object is a generalized model of the
discrete memoryless one.  Let  ${\cal X}$ be a finite set of values
of RV  $X$, and  ${\cal S}$ is an alphabet of states of the object.

  $M$ possible conditional PDs   of the characteristic $X$ of the
object  depending on values $s$ of  states,    are given:
$$
W_m\df \{W_m(x|s), \,\,\,x\in {\cal X},\,\,\,\, s\in {\cal S}
\},\,\,\,\,m=\overline {1,M},... |{\cal S}\geq 1|,
$$
but it is not known which of these
 alternative hypotheses
$W_m$, $m=\overline {1,M}$,  is   real PD of the object.  The
statistician must select one among $M$ hypotheses, or he can
withdraw any judgement. It is possible for instance when it is
supposed that real PD is not in  the family of $M$ given PDs.
 An answer must be given   using the vector of
results of $N$ independent experiments ${\bf x}\df(x_1,x_2,...x_N)$
and the  vector of states of the object ${\bf s}\df
(s_1,s_2,...,s_N)$, $s_n\in {\cal S}$, $n=\overline {1,N}$.

 The procedure of decision making is a non-randomized test $\varphi_N({\bf x}, {\bf s})$,
 it can be defined  by division of
 the sample space ${\cal X}^N$   for each $\bf s$ on $M+1$ disjoint subsets
   ${\cal A}^N_m({\bf s})=\{{\bf x:}\,\,\, \varphi_N({\bf x}, {\bf s})=m\}$, $m=\overline {1,M+1}$.
   The  set  ${\cal A}^N_l({\bf s})$, $l=\overline {1,M}$, consists of  vectors $\bf x$
   for which the hypothesis $W_l$
    is adopted, and ${\cal A}^N_{M+1}({\bf s})$ includes   vectors for which the statistician
    refuses
    to take a certain answer.

 We study   the probabilities of the erroneous acceptance
  of hypothesis $W_l$ provided that $W_m$ is true
\begin{eqnarray}
\label {54} \alpha_{l|m}(\varphi_N)\df \max\limits_{{\bf s}\in {\cal
S}^N}W^N_m\left ({\cal A}_l^{N}({\bf s})|{\bf s}\right ),\,\,
m,l=\overline{1,M},\,\,\,m\not= l.
\end{eqnarray}
 When
decision is declined, but  hypothesis $W_m$ is true, we consider the
following probability of error:
$$
\alpha_{M+1|m}(\varphi_N)\df \max\limits_{{\bf s}\in {\cal
S}^N}W^N_m\left ({\cal A}_{M+1}^N({\bf s})|{\bf s}\right ).
$$

If the hypothesis $W_m$ is true, but it is not accepted, or
equivalently   while the statistician accepted one of hypotheses
$W_l$, $l=\overline {1, M}$, $l\not= m$, or refused to make
decision, then the  probability of error is the following:
\begin{eqnarray}
\label {55} \alpha_{m|m}(\varphi_N)\df \sum\limits_{l:\,l\not=
m}\alpha_{l|m}(\varphi_N)=   \max\limits_{{\bf s}\in {\cal
S}^N}W^N_m\left ({\overline {{\cal A}^N_m({\bf s})}}|{\bf s}\right
), \,\,\,\, m=\overline {1,M}. \end{eqnarray}

Corresponding reliabilities  are defined similarly by to  (\ref
{5}):
\begin{eqnarray}
\label {56} E_{l|m}(\varphi)\df\overline {\lim\limits_{N\to
\infty}}\left \{-\frac{1}{N}\log\alpha_{l|m}(\varphi_N)\right \},
\,\, m=\overline {1,M},\,\,\,\,\,l=\overline {1,M+1}.
\end {eqnarray}

It also follows   that for every test $\varphi$
\begin{eqnarray}
\label {57} E_{m|m}(\varphi)={\mathop {\min}\limits_{l=\overline
{1,M+1},\,\, l\not= m}}E_{l|m}(\varphi),\,\,\,m=\overline {1,M}.
\end
{eqnarray} The matrix
$$
{\bf E}(\varphi)=\left (
\begin{array}{c}
E_{1| 1}\,\, \ldots \, E_{ l|1}\, \ldots \,\, E_{M|1 },\,\,E_{M+1|1 }\\
\ldots \ldots \ldots \ldots  \ldots \ldots \ldots \\
E_{ 1|m} \ldots E_{ l|m} \ldots   E_{ M|m},\,\,E_{ M+1|m}\\
\ldots \ldots \ldots \ldots  \ldots \ldots \ldots \\
E_{1|M  } \ldots  E_{ l|M } \ldots  E_{M| M}\,\,E_{ M+1|M}
\end{array}
\right )
$$
 is  the reliabilities matrix of the tests sequence $\varphi$
 for the described model.

 We call the  test LAO for this model  if for given positive values of certain  $M$
 elements of the matrix
${\bf E}(\varphi)$ the procedure  provides maximal values for other
elements of it.

 For construction of    LAO  test
 positive elements $E_{1|1}$, $...$, $E_{M|M}$ are supposed to be  given
preliminarily. The optimal dependence of
 error exponents was determined in   \cite {EHA}.
 This result     can be easily generalized for  the case of an arbitrarily varying Markov source.

%
%
%

\subsection { Multiple Hypotheses  LAO Testing With Rejection of Decision for Many  Independent Objects}
For brevity we consider the problem for two objects,  the
generalization of the problem  for $K$
  independent objects will be discussed   along  the text.

  Let $X_1$ and $X_2$ be independent RVs taking values in the same finite set ${\cal X}$
with one of $M$ PDs $G_{m}\in {\cal P}({\cal X})$, $m=\overline
{1,M}$. These RVs are the characteristics of the corresponding
independent objects. The random vector $(X_1,X_2)$ assumes values
$(x^1,x^2)\in {\cal X}\times {\cal X}$. Let $({\bf x^1},{\bf
x^2})\df $ $\left
((x_{1}^{1},x_{1}^{2}),...,(x_{n}^{1},x_{n}^{2}),...,
(x_{N}^{1},x_{N}^{2})\right )$, $x_n^k\in {\cal X},$ $k=\overline
{1,2},
 n=\overline {1,N}$,  be a vector of results of $N$ independent observations of the pair of RVs   $(X_1,X_2)$.
 On the base of   observed data the test has  to determine   unknown PDs of the objects or  withdraw
any judgement.
  The selection for each object should  be
 made from   the same set of hypotheses: $G_m$, $m=\overline {1,M}$.
We call  this  procedure  the  compound test for two  objects and
denote it by $\Phi_N$, it can be composed of   two  individual tests
$ \varphi^{1}_N$, $\varphi^{2}_N$   for corresponding  objects. The
test $\varphi^{i}_N$, $i=\overline {1,2}$, can be defined by
division of the space  ${\cal X}^N$ into $M+1$ disjoint subsets
${\cal A}_m^i$, $m=\overline {1, M+1}$. The set ${\cal A}_m^i$,
$m=\overline {1,M}$ contains all vectors ${\bf x}^i$ for which the
hypothesis $G_m$ is adopted and ${\cal A}_{M+1}^i$ includes all
vectors for which the test refuses to take a certain answer. Hence
$\Phi_N$ is division of the space ${\cal X}^N\times {\cal X}^N$ into
$(M+1)^2$ subsets ${\cal A}_{m_1,m_2}={\cal A}_{m_{1}}^1\times {\cal
A}_{m_2}^2$, $m_i=\overline {1,M+1}$. We again  denote the infinite
sequences of  tests by $ {\Phi} $, $\varphi^{1}$, $\varphi^{2}$.

Let   $\alpha_{l_1,l_2|m_1,m_2}(\Phi_N)$ be the probability of the
erroneous acceptance    of  the pair of hypotheses
$(G_{l_{1}},G_{l_{2}})$ by the test $\Phi_N$ provided that the pair
of hypotheses $(G_{m_{1}},G_{m_{2}})$ is  true, where
$(m_1,m_2)\not=(l_1,l_2)$, $m_i= \overline {1,M}$, $l_i=\overline
{1,M}$, $i=\overline {1,2}$:
\begin {align}
\alpha_{l_1,l_2|m_1,m_2}(\Phi_N)&=G_{m_1}\circ G_{m_2}\left ({\cal
A}_{l_1,l_2}\right )\nonumber\\
 &=G_{m_1}^N({\cal A}_{l_1})\cdot
G_{m_2}^N({\cal A}_{l_2}).\nonumber
\end {align}

When   the pair of  hypotheses  $(G_{m_{1}},G_{m_{2}})$, $m_1,
m_2=\overline {1,M}$ is true, but we decline the decision the
corresponding probabilities of errors are:
\begin {align}
\alpha_{M+1,M+1|m_1,m_2}(\Phi_N)&=G_{m_1}\circ G_{m_2}({\cal
A}_{M+1,M+1})\nonumber\\ & =G_{m_1}^N({\cal A}_{M+1}^1)\cdot
G_{m_2}^N({\cal A}_{M+1}^2).\nonumber
 \end  {align}
 or
$$ \alpha_{ M+1,l_2|m_1,m_2}(\Phi_N)=G_{m_1}^N({\cal
A}_{M+1}^1)\cdot G_{m_2}^N({\cal A}_{l_2}^2)$$
 or
$$ \alpha_{l_1,M+1|m_1,m_2}(\Phi_N)=G_{m_1}^N({\cal A}_{l_1}^1)\cdot
G_{m_2}^N({\cal A}_{M+1}^2).
$$

 If the  pair of hypotheses
$(G_{m_{1}},G_{m_{2}})$ is true, but it is not accepted, or
equivalently   while the statistician accepted one of hypotheses
$(G_{l_{1}},G_{l_{2}})$,  or refused to make decision, then the
probability of error is the following:
\begin {align}
\label{60}\alpha_{m_1,m_2|m_1,m_2}(\Phi_N)&=
\sum\limits_{(l_1,l_2)\not=(m_1,m_2)}
 \alpha_{l_1,l_2|m_1,m_2}(\Phi_N),\\
&\,\,l_i=\overline {1, M+1}, \,\,m_i=\overline {1,M},\,\,i=\overline
{1,2}.\nonumber
\end {align}

We  study    reliabilities
  $E_{l_1,l_2|m_1,m_2}(\Phi)$ of the sequence of
tests $\Phi$,
 \begin {align}
E_{l_1,l_2|m_1,m_2}(\Phi)&\df
 \label {61} \overline{\lim\limits_{N\to\infty}}-\frac{1}{N}
\log\alpha_{l_1,l_2|m_1,m_2} (\Phi_N),\\
 &\,\,\,m_i,= \overline {1,M},\,\,l_i=\overline
{1,M+1},\,\,\,\,\,i=\overline {1,2}. \nonumber
\end {align}
 Definitions (\ref {60}) and (\ref {61})  imply   that
\begin {align}
\label {62}
E_{m_1,m_2|m_1,m_2}(\Phi)&=\min\limits_{(l_1,l_2)\not=(m_1,m_2)}
E_{l_1,l_2|m_1,m_2}(\Phi).
\end
{align}

We can  erect the LAO test from the set of compound tests when $2M$
strictly positive elements of the reliability matrix $E_{M+1,m|m,m}$
and $E_{m,M+1|m,m}$, $m=\overline {1,M}$, are preliminarily given
(see \cite {EHA}).


{\bf Remark 6.1:}{ It is necessary to note that the problem of
reliabilities investigation for LAO  testing of  many hypotheses
with possibility of  rejection of decision for the model consisting
of two or more independent objects can not be  solved  by the direct
method of renumbering.}

\section  {Conclusion and Open Problems}

``A broad class of statistical problems arises in the framework of
hypothesis testing in the spirit of
  identification for different kinds of sources, with complete or
  partial  side information or without it. Paper \cite {R} is a
  start." \cite {Rnew}.

 In this paper, we exposed solutions of a part of possible problems concerning  algorithms of
  distributions optimal testing
 for certain classes of one, or multiple objects. For the
 same models PD optimal identification  is discussed again in the
 spirit of error probability exponents optimal dependence. But these
 investigations can be continued in plenty directions.

 Some  problems   formulated in \cite {R} and \cite {HMA},
 particularly, concerning the remote statistical inference
 formulated by Berger \cite {ber}, examined  in part by Ahlwede and Csisz\'ar
 \cite {AC} and Han and Amari \cite {HnA} still rest open.

 All our results concern with discrete distribution, it is
necessary to study many objects with general distributions as in
\cite {13}. For multiple  objects    multistage and sequential
testing \cite {DB} can be also considered. Problems for many objects
are
  present in statistics with fuzzy data \cite {24}, bayessian detection of multiple hypotheses testing \cite {LJ}
   and geometric interpretations of tests \cite {Mbw}.

\section {Appendix}

 {\it Proof of Theorem 3.1:} Probability $G_m^N({\bf x})$ for   ${\bf x}\in {\cal
T}^N_{Q}(X)$
   can be presented as follows:
\begin {align}
G_m^N({\bf x}) &= \prod\limits_{n=1}^N G_m(x_n)\nonumber\\
&= \prod\limits_{x }G_m(x )^{N(x \mid {\bf x})}\nonumber\\
&= \prod\limits_{x }G_m(x)^{NQ(x))} \nonumber\\
 &=  \exp\left \{N \sum\limits_{x}\left (-Q(x)
 \log\frac{Q(x)}{G_m (x)}+Q(x)\log{Q(x)}\right )\right \}  \nonumber\\
 &= \exp{\left \{-N[D(Q \parallel G_m )+H_{Q}(X)]\right
\}}.\label {70}
\end {align}
 Let us consider the  sequence of tests
$\varphi^*_N({\bf x})$ defined  by the  sets
\begin {align}
\label {71} {\cal B}_m^{(N)}&\df \bigcup_{P\in {\cal R}^{(N)}_{m}}
{\cal T}^N_{Q}(X),\,\,\,m=\overline {1,M}.
\end {align}
Each $\bf x$ is in one and only in one of ${\cal B}_m^{(N)}$, that
is
$$
{\cal B}_l^{(N)}\bigcap {\cal B}_m^{(N)}=\emptyset,\,\,\,l\not= m,
\,\,\,\,\,\,\, \mbox {and}\,\,\,\,\,\,\,\,\,\bigcup_{m=1}^M{\cal
B}^{(N)}_m={\cal X}^N.
$$
Really, for $l=\overline {1,M-2}$,  $m=\overline {2,M-1}$,  for each
$l<m$   let us  consider arbitrary ${\bf x}\in {\cal B}_l^{(N)}$. It
follows  from (\ref {8}) and (\ref {10}) that there exists type
$Q\in {\cal Q}^N({\cal X})$ such that $D(Q||G_l)\leq E_{l|l}$ and
${\bf x}\in {\cal T}^N_{Q}(X)$. From (\ref {15}) and (\ref {12}) we
have $E_{m|m}<E_{l|m}^*(E_{l|l})<D(Q||G_m)$. From definition of
${\cal B}_m^{(N)}$ we see that ${\bf x}\notin {\cal B}_m^{(N)}$.
Definitions
 (\ref {13}), (\ref {15}) and (\ref {10})  show also that
$$
{\cal B}_M^{(N)} \bigcap {\cal B}_m^{(N)}
=\emptyset,\,\,\,m=\overline {1,M-1}.
$$
Now, let us  remark   that    for $m=\overline {1,M-1}$, using (\ref
{1}),  (\ref {2}),  (\ref {4})--(\ref {6}) and (\ref {70}) we can
estimate $\alpha_{m|m}^{(N)}(\varphi^*)$ as follows:

\begin {align}
\alpha_{m|m} (\varphi^*_N)  &=   G_m^N\left (\overline {{\cal B}_m^{(N)}}\right )\nonumber\\
&= G_m^N\left (\bigcup_{Q:D(Q ||G_m)> E_{m|m}}{\cal T}^N_{Q}(X)
\right )\nonumber\\
&\leq  (N+1)^{|{\cal X}|} \sup\limits_{Q :D(Q||G_m)>
E_{m|m}}G_m({\cal T}^N_{Q}(X))\nonumber\\
&\leq  (N+1)^{|{\cal X}|}\sup\limits_{Q:D(Q||G_m)>
E_{m|m}}\exp\{-ND(Q||G_m)\}\nonumber\\
&\leq \exp \left \{ -N[ \inf\limits_{Q :D(Q||G_m) >
E_{m|m}}D(Q||G_m)-o_N(1)]\} \right
\}\nonumber\\
& \leq  \exp \left \{-N[E_{m|m}-o_N(1)]\right \},\nonumber\\
\nonumber
\end {align}
where $o_N(1)\to 0$ with $N\to \infty$.

 For $l=\overline {1,M-1}$,
$m=\overline {1,M}$, $l\not=m$, using (\ref {1}),  (\ref {2}), (\ref
{4})--(\ref {6}) and (\ref {70}),  we can obtain  similar estimates:
\begin {align}
\label {72} \alpha_{l|m} (\varphi^*_N)  &=  G_m^N\left ({\cal
B}_l^{(N)}\right )\nonumber\\
&=  G_m^N\left (\bigcup_{Q\,:D(Q||G_l)  E_{l|l}}{\cal T}_{Q}^N   (X) \right)\nonumber\\
&\leq  (N+1)^{|{\cal X}|}\sup\limits_{Q:D(Q||G_m)\leq
E_{l|l}}G_m^N({\cal T}_{Q}^N(X)\nonumber\\
 &\leq (N+1)^{|{\cal X}|}\sup\limits_{Q:D(Q||G_m)\leq E_{l|l}}\exp
\{-ND(Q||G_m)\}\nonumber\\
&= \exp\left \{-N\left (\inf\limits_{Q:D(Q||G_m)\leq
E_{l|l}}D(Q||G_m) - o_N(1)\right )\right \}.\\
\nonumber
\end {align}
Now let us prove the inverse inequality:
\begin {align}
\alpha_{l|m} (\varphi^*_N) & =   G_m^N\left ({\cal B}_l^{(N)}\right )\nonumber\\
&= G_m^N\left (\bigcup_{Q : D(Q||G_l)\leq
E_{l|l}}{\cal T}_{Q}^N(X)\right )\nonumber\\
& \geq  \sup\limits_{Q :
D(Q||G_l)\leq E_{l|l}}G^N_m({\cal T}_{Q}(X)\nonumber\\
\label {73}
 & \geq (N+1)^{-|{\cal X}|}\sup\limits_{Q : D(Q||G_l) \leq
E_{l|l}}\exp \left \{-ND\left (Q||G_m\right )\right \}\nonumber\\
&= \exp\left \{-N\left (\inf\limits_{Q: D(Q||G_l)\leq
E_{l|l}}D(Q||G_m)+o_N(1)\right )\right \}.
\end {align}
According to the definition (\ref {6}) $E_{l|m}(\varphi^*)=\overline
{\lim\limits_{N\to \infty}}\left \{-N^{-1}\log\alpha_{l|m}
(\varphi^*_N)\right \}$, taking into account (\ref {72}), (\ref
{73})
and the continuity of the functional $D(Q||G_l)$ we obtain that\\
$\lim\limits_{N\to \infty}\left
\{-N^{-1}\log\alpha_{m|l}(\varphi^*_N)\right \}$ exists and in
correspondence with (\ref {12}) equals to $E_{m|l}^*$. Thus
$E_{l|m}(\varphi^*)=E_{l|m}^*$, $m=\overline {1,M}$, $l=\overline
{1, M-1}$, $l\not= m$. Similarly we can obtain upper and lower
bounds for $\alpha_{M|m} (\varphi^*_N)$, $m=\overline {1,M}$.
Applying the same resonnement we get that  the reliability
$E_{M|m}(\varphi^*)=E^*_{M|m}$.
 By the definition  (\ref {7}) $E_{M|M}(\varphi^*)=E^*_{M|M}$.
  The proof of the first part of the theorem will be accomplished if
we demonstrate that the sequence of  tests $\varphi^*$ is LAO, that
is for given  $E_{1|1},...,E_{M-1|M-1}$  and  every  sequence of
tests ${\varphi}$   for all $l,m\in\overline {1,M}$,
$E_{m|l}(\varphi)\leq E^*_{m|l}.$

 Let us consider any other sequence $\varphi^{**}$ of tests which are defined
  by the sets ${\cal D}_1^ {(N)},...,{\cal
D}_M^{(N)}$ such that
\begin {equation}
\label {M} E_{l|m}(\varphi^{**})\geq E_{l|m}^*, \,\,\,m,l=\overline
{1,M}.
\end{equation}

 These conditions are equivalent for $N$ large enough to the inequalities
\begin {equation}
\label {74} \alpha_{l|m}(\varphi^{**}_N)\leq
\alpha_{l|m}(\varphi^{*}_N), \,\, m,l=\overline {1,M}.
\end {equation}
Let us examine the sets $D_m^{(N)}\bigcap  {\cal B}_m^{(N)}$,
$m=\overline {1, M}$. This intersection cannot be empty, because in
that case
\begin {align}
\alpha_{m|m}(\varphi^{**}_N)& = G_m^N\left ({\overline {\cal D}}_m^{(N)}\right )\nonumber\\
&\geq  G_m^N\left ({\cal B}_m^{(N)}\right )\nonumber\\
&\geq  \max\limits _{Q
:D(Q||G_m)\leq E_{m|m} }G_m^N\left ({\cal T}_{Q}^N(X)\right )\nonumber\\
&\geq \exp\{-N(E_{m|m}+o_N(1))\}.\nonumber
\end {align}
Let us show that ${\cal D}_l^{(N)}\bigcap {\cal
B}_m^{(N)}=\emptyset$, $m,l= \overline {1, M-1}$, $m\not=l$. If
there  exists $Q$ such that  $D(Q||G_m)\leq E_{m|m}$ and ${\cal
T}_{Q}^N(X)\in {\cal D}_l^{(N)}$, then
\begin {align}
\alpha_{l|m}(\varphi^{**}_N)&= G_m^N\left ({\cal D}_l^{(N)}\right )\nonumber\\
&>  G_m^N\left ({\cal
T}_{Q}^N(X)\right )\nonumber\\
&\geq  \exp \{-N(E_{m|m} +o_N (1))\}.\nonumber
\end {align}
When $\emptyset\not= {\cal D}_l^{(N)} \bigcap {\cal T}_{Q}^N(X)\not=
{\cal T}_{Q}^N(X)$, we also obtain that
 \begin {eqnarray*}
\alpha_{l|m}(\varphi^{**}_N)& = &G_m^N({\cal D}_l^{(N)}\nonumber\\
 & > &G_m^N({\cal D}_l^{(N)}
 \bigcap {\cal T}_{Q}^N(X)\nonumber\\
 &\geq&\exp \{-N(E_{m|m} +o_N (1))\}.
 \end {eqnarray*}
Thus it follows  that $E_{l|m}(\varphi^{**})\leq  E_{m|m}$, which in
turn according to (\ref {7}) provides that $E_{l|m}(\varphi^{**})=
E_{m|m}$. From condition (\ref {15}) it follows that $E_{m|m}<
E_{l|m}^*$ , for all  $l=\overline {1, m-1}$, hence
$E_{l|m}(\varphi^{**})< E_{l|m}^*$ for all  $l=\overline {1, m-1}$,
which contradicts to (\ref {M}). Hence we obtain that ${\cal
D}_m^{(N)}\bigcap  {\cal B}_m^{(N)}={\cal B}_m^{(N)}$ for
$m=\overline {1, M-1}$. The intersection $D_m^{(N)}\bigcap {\cal
B}_M^{(N)}$ is empty too, because otherwise
$$
\alpha_{M|m}(\varphi^{**}_N)\geq \alpha_{M|m}(\varphi^{*}_N),
$$
which contradicts to (\ref {74}), hence  ${\cal
D}_m^{(N)}=B_m^{(N)}$, $m=\overline {1,M}$.

The proof of the second part  of the Theorem 3.1 is simple. If one
of the conditions (\ref {15}) is violated, then from  (\ref
{12})--(\ref {14}) it follows that at least one of the elements
$E_{m|l}$ is equal to 0. For example, let   in (\ref {15}) the
$m$-th condition be  violated. It  means that
$E_{m|m}\geq\min\limits_{l=\overline {m+1,M}}D(G_l||G_m)$, then
there exists $l^*\in \overline {m+1,M}$ such that $E_{m|m}\geq
D(G_l^*||G_m).$ From latter and (\ref {12}) we obtain that
$E^*_{m|l}=0$.

The theorem is proved.
 \vspace {3mm}

 {\it Proof of Lemma 4.1:}  It follows
from the independence of the objects  that
\begin{equation}
\hspace {-2,1cm}\alpha_{
l_1,l_2,l_3|m_1,m_2,m_3}(\Phi_N)=\prod\limits_{
i=1}^3\alpha_{l_i|m_i}(\varphi_{N}^i), \,\,\,\mbox{if}\,\,\,\,
m_i\not=l_i,
 \label{13.a}
\end{equation}
$$
\alpha_{ l_1,l_2,l_3|m_1,m_2,m_3}(\Phi_N)=\left (1-\alpha_{l_k|m_k}
(\varphi_{N}^k)\right ) \prod\limits_{ i\in [[1,2,3]-k]}
\alpha_{l_i|m_i}(\varphi_{N}^i),
$$
\begin{equation}
 m_{k}=l_{k},\,\, m_i\not=l_i, \,\, k=\overline {1,3},\,\,\,i\not=
 k,
 \label{13.b}
\end{equation}
$$
\alpha_{ l_1,l_2,l_3|m_1,m_2,m_3}(\Phi_N)=
\alpha_{l_i|m_i}(\varphi_{N}^i)\prod\limits_{ k\in [[1,2,3]-i]}
\left (1-\alpha_{l_k|m_k}(\varphi_{N}^i)\right ),
$$
\begin{equation}
m_{k}=l_{k},\,\, m_i\not=l_i,\,\,i=\overline {1,3}.
 \label{13.c}
\end{equation}

Remark that here we    consider also  the  probabilities of right
(not erroneous) decisions. Because $E_{l|m}(\varphi^i)$ are strictly
positive  then  the  error probability\\
$\alpha_{l|m}(\varphi^i_N)$ tends to zero, when $N\longrightarrow
\infty$. According this fact we have
\begin{eqnarray}
\overline{\lim\limits_{N\to \infty}}\left \{-\frac{1}{N}\log\left
(1-\alpha_{l|m}(\varphi^i_N)\right )\right \}&=&
\overline{\lim\limits_{N\to \infty}}\frac {\alpha_{l|m}\left
(\varphi^i_N\right )}{N}\frac{\log\left
(1-\alpha_{l|m}(\varphi^i_N)\right
)}{-\alpha_{l|m}(\varphi^i_N)}\nonumber\\
\label{14.a} &=&0.\\ \nonumber
\end{eqnarray}
From  definitions (\ref {18}), equalities  (\ref {13.a})--(\ref
{13.c}), applying  (\ref {14.a}) we obtain
 relations  (\ref {20})--(\ref{22}).

The Lemma is proved.

 \vspace {3mm}

{\it Proof of Theorem 4.1:} The test $\Phi^*=(\varphi^{1,*},
\varphi^{2,*}, \varphi^{3,*})$, where $\varphi^{i,*}$, $i=\overline
{1,3}$ are LAO tests of objects $X_i$, belongs to the set ${\cal
D}$. Our aim is to prove  that such $\Phi^*$ is a compound   LAO
test.
  Conditions  (\ref {28})--(\ref {31})   imply  that inequalities
analogous to (\ref {15}) hold simultaneously for  tests  for three
separate objects.

Let   the test  $\Phi\in {\cal D}$ be  such  that
$E_{M,M,M|m,M,M}(\Phi)=E_{M,M,M|m,M,M}$\\
$E_{M,M,M|M,m,M}(\Phi)=E_{M,M,M|M,m,M}$,  \,and
$E_{M,M,M|m,M,M}(\Phi)=E_{M,M,M|M,M,m}$, $m= \overline {1,M-1}$.

Taking into account  (\ref {23_1})--(\ref {23_3}) we can see that
conditions (\ref {28})--(\ref {31}) for every $m=\overline {1,M-1}$
may be replaced by  the following inequalities:
\begin {equation}
E_{M|m}(\varphi^i)<\min \left [\min\limits_{l=\overline
{1,m-1}}\inf\limits_{Q:D(Q||G_m)\leq
E_{M|m}(\varphi^i)}D(Q||G_l),\min\limits_{l=\overline
{m+1,M}}D(G_l||G_m)\right ]. \label{123}
\end {equation}

According to Remark 3.1 for  LAO test $\varphi ^{i,*}$, $i=\overline
{1,3}$, we obtain that (\ref {123})  meets  conditions
 (\ref {15}) of  Theorem 3.1 for
 each  test  $\Phi\in {\cal D}$,
$E_{m|m}(\varphi^i)>0$, $i=\overline {1,3}$, hence it follows from
(\ref {6}) that $E_{m|l}(\varphi^i) $ are also strictly positive.
Thus for a test $\Phi\in {\cal D}$ conditions of Lemma 4.1 are
fulfilled and the elements of the reliability matrix ${\bf E}(\Phi)$
coincide with elements of matrix ${\bf E}(\varphi^i)$, $i=\overline
{1,3}$, or sums of them. Then from definition of LAO test it follows
that $E_{l|m }(\varphi^i)\leq E_{l|m } (\varphi ^{i,*})$, then
$E_{l_1,l_2,l_3|m_1,m_2,m_3} (\Phi) \leq E_{l_1,l_2,l_3|m_1,m_2,m_3}
(\Phi^*) $. Consequently $\Phi^*$ is a LAO test and
$E_{l_1,l_2,l_3|m_1,m_2,m_3} (\Phi^*) $ verify $(\ref {24})$--$(\ref
  {27})$.

 b) When even one of the inequalities $(\ref {28})$--$(\ref {31})$
 is violated, then at least one of inequalities (\ref {123}) is  violated. Then from
 Theorem 3.2 one of elements $E_{m|l}(\varphi^{i,*})$ is equal to
 zero. Suppose $E_{3|2}(\varphi^{1,*})=0$, then the elements\\
 $E_{3,m,l|2,m,l}(\Phi^*)=E_{3|2}(\varphi^{1,*})=0$.

 The Theorem   is proved.

\end{document}